%% file: paper_revised3.tex
\documentclass[article,3p]{elsarticle} 

\usepackage{amsmath,amsfonts,amssymb,amsbsy,amscd}

\usepackage{mathrsfs,array}
\usepackage{verbatim}
\usepackage{lipsum}
\usepackage{subcaption}
\usepackage{amsthm}
\usepackage{tikz}

\newdefinition{rmk}{Remark}
\newproof{pf}{Proof}
\newproof{pot}{Proof of Theorem \ref{thm2}}

\input{declarations}

\biboptions{numbers,sort&compress}

\usepackage{hyperref}
\hypersetup{
    colorlinks,
    citecolor=ck,
    filecolor=ck,
    linkcolor=ck,
    urlcolor=ck
}

\newcommand{\Marco}[1] { {\color{blue}  #1} }
\newcommand{\Ido}[1] { {\color{violet}  #1} }
\newcommand{\Derive}[1] { {\color{green}  #1} }
\newcommand{\TBD}[1] { {\color{red} #1} }

\renewcommand{\Marco}[1] {}
\renewcommand{\Ido}[1] {}
\renewcommand{\Derive}[1] {}
\renewcommand{\TBD}[1] {}

\usepackage[final]{changes}

\definechangesauthor[]{Rev.1}
\definechangesauthor[]{Rev.2}
\definechangesauthor[]{Both Rev}
\definechangesauthor[]{Authors}

\definechangesauthor[color=black]{Rev2.1}
\definechangesauthor[color=black]{Rev2.2}
\definechangesauthor[color=orange]{Rev3.2}

\newcommand{\opt}{bw}
\renewcommand{\opt}{color}
\renewcommand{\opt}{color_dash}
\journal{}
\begin{document}

\makeatletter
\def\ps@pprintTitle{%
   \let\@oddhead\@empty
   \let\@evenhead\@empty
   \let\@oddfoot\@empty
   \let\@evenfoot\@oddfoot
}
\makeatother
\begin{frontmatter}

        
\title{Toward free-surface flow simulations with correct energy evolution:\\ an isogeometric level-set approach with monolithic time-integration}


\author{I. Akkerman\corref{cor1}}
\ead{i.akkerman@tudelft.nl}
\cortext[cor1]{Corresponding author}
\author{M.F.P. ten Eikelder}
\ead{m.f.p.teneikelder@tudelft.nl}

\address{Delft University of Technology, \\
Department of Mechanical, Maritime and Materials Engineering,\\
 P.O. Box 5, 2600 AA Delft, The Netherlands}
\date{}

\begin{abstract}
This paper presents a new monolithic free-surface formulation that exhibits correct kinetic and potential energy behavior.
\added[id=Rev2.1]{We focus in particular on the temporal energy behavior of two-fluids flow with varying densities.}
\added[id=Rev3.2]{Correct energy behavior here means that the actual energy evolution of the numerical solution matches the 
evolution as predicted by the discrete two--fluid equations.}
We adopt the level-set method to describe the two-fluid surface.
To ensure the correct energy behavior we augment the interface convection equation with kinetic and potential energy constraints.
We solve the resulting formulation consisting of the fluid and interface equations in a monolithic fashion using a recently proposed level-set method 
[I. Akkerman. Monotone level-sets on arbitrary meshes without redistancing. Computers $\&$ Fluids, 146:74-85, 2017.]. 
For the spatial discretization divergence-conforming NURBS are adopted.
The resulting discrete equations are solved with a quasi-newton method which partially decouples the constraints from the rest of the problem.

\added[id=Both Rev]{As we focus on the energy behavior of  time integration in case of  varying densities, 
we restrict ourselves to low-Reynolds-number flow allowing simple Galerkin discretizations. 
High-Reynolds-number two-fluid flows that require stabilization are beyond the scope of the current paper.}
The simulation of a dambreak problem numerically supports the correct energy behavior of the proposed methodology.
The proposed methodology improves the solution quality significantly upon a more traditional approach.
Due to the excellent accuracy per degree of freedom one can suffice with a much lower resolution.
\end{abstract}


\end{frontmatter}

\section{Introduction}
\label{sec:intro}
Free-surface problems are ubiquitous in science and engineering, in particular problems involving
air-water surfaces are often encountered in a maritime, offshore or coastal engineering.
Numerical methods developed to simulate free-surface problems often use the level-set method, originally proposed in \cite{Osher_Sethian_88}, to describe the evolution of the free-surface.
The level-set approach avoids the possibility of negative densities which could break down the entire computation.
Variational formulations of convection problems generally (excluding some notable exceptions \cite{MIZUKAMI1985, BMPR96, EVANS2009}) do not satisfy the maximum principle.
This means that under- and overshoots of the density profile can appear, especially when dealing with large density jumps (e.g. air-water flows).
Since the level-set method precludes this, it is a very popular method for all sorts of evolving interface problems, see for instance the review papers \cite{OSHER2001,GIBOU2018,SUSSMAN1994,AMW98,ScZa99}.
The methodology is particularly often employed when a finite element method is used for the discretization \cite{Akin05a,LoYaOn06,ElCo07,LERC10,CrCeTe07,CrCeTe07b,KAFB10}.
Also the combination of isogeometric analysis \cite{CoHuBa09} with level-sets has been explored \cite{ABKF11}.
Moreover, these methodologies have also been used to perform fluid-structure interaction (FSI) in conjunction with a free-surface, see e.g. \cite{WALHORN2005, AkBaBeFaKe12,ADKSB12}. 

The combination of a fluid and a structure is often solved monolithically in the FSI community.
This approach has clear advantages in terms of stability which translates into improved robustness and efficiency.
The advantage of monolithic coupling is described in a general setting in \cite{Felippa01}, 
while the gain of monolithic free-surface modeling is demonstrated in \cite{AkBaBeFaKe12}.
In the latter paper the free-surface/rigid body problem is formulated and solved in a strongly coupled way.
Comparing with a staggered Navier-Stokes/level-set convection formulation, this significantly improves the stability of the methodology.
The staggered formulation was shown to create artificial energy, which rendered the solution completely useless at some point in time. 
For ease of implementation, redistancing of the level-set\footnote{The level-set methodology requires a scaling and redistancing step of the level-set to have control over the size of the smoothing region around the interface.} 
and a mass correction step were kept out of the main iteration loop.
Instead, these corrections to the interface were applied after the main solve.
This leads to small errors in the conserved properties momentum and energy.

The creation of artificial energy in the numerical formulation is evidently considered unfavorable. 
The numerical energy plays a fundamental role in the numerical stability of the method, which was already recognized in \cite{lax1974}.
Despite that fact that artificial energy creation could lead to useless solutions, many methods developed for the simulation of two-fluid flows can unfortunately create artificial energy.
Moreover, even in the mono-fluid case the popular stabilized finite element methods can create artificial energy \cite{EiAk17i, EiAk17ii}.
These papers correct this imperfection in the single-fluid case.

In this paper we develop a free-surface formulation based on level-sets that has guaranteed correct energy behavior.
That is the conservation of energy in the inviscid case, and guaranteed physical energy decay when viscosity is present.
We solve the fluid and interface evolution problem, including redistancing and mass correction, in a monolithic fashion.
Therefore we make use of the novel level-set formulation with an efficient and robust redistancing approach developed in \cite{AKKERMAN2017}.
The standard discretization appears to have a mismatch between discretized and continuous kinetic and potential energy.
To correct this mismatch and to ensure the correct behavior of these energies we augment the convection equation with the required constraints.
Demanding mass conservation results in an additional constraint.
To enforce these constraints we utilize the method of Lagrange-multipliers.

The goal of the current paper is to develop a time-integration procedure for two-fluid flow that exhibits the correct energy behavior. 
This procedure should not artificially destroy or create energy at the interface due to changing densities. 
In order to keep the focus on this, we do not incorporate stabilization techniques. 
Hence, this paper deals with low Reynolds number flows.
A follow-up paper addresses the high Reynolds number case.

The outline of the remainder of this paper is as follows.
Section \ref {sec:continuous} presents the continuous form of the governing fluid equations, both in strong and weak form.
Furthermore, it discusses the conservation of mass, momentum and energy in the continuous setting.
Section \ref{sec:disc_cons_form} provides the standard discretization. 
It closely mimics the weak continuous form in a discrete setting.
The standard formulation is shown not to have the correct interface evolution 
in order to guarantee correct evolution of mass, kinetic energy and potential energy.
Moreover, this section also presents the employed level-set method. 
The methodology is presented in the isogeometric analysis framework which is shown to be beneficial for the energy evolution behavior.
Next, section \ref{sec:EC_form} presents our novel constrained method that corrects these discrepancies.
In contrast to the standard discretization, this approach displays the correct mass, potential energy and kinetic energy behavior.
The discrete formulation is solved using a quasi-newton solver.
The numerical verification on a prototype dam-break problem is presented in section \ref{sec:results}.
It compares the energy evolution of the standard method and the newly proposed method with a benchmark convective approach.
In the final section \ref{sec:conclusion} we draw conclusions and discuss some avenues for future work.

\section{Continuous formulation}
\label{sec:continuous}
\Marco{In principle the organization of this section is good.}

\subsection{The governing equations}
Let $\Omega \subset \mathbb{R}^d$, $d=2,3$, denote the spatial domain with boundary $\partial \Omega = \Gamma$. 
The problem under consideration consists of solving the incompressible Navier-Stokes equations dictating the two-fluid flow:
\begin{subequations}\label{eq:strong_cons_ns}
\begin{alignat}{3}
 \partial_t (\rho \bu) + \nabla \cdot (\rho \bu\otimes \bu) 
 + \nabla  p - \nabla \cdot 2\mu \nabla^s \bu &=  \rho \mathbf{g} &\quad \text{in}& \quad \Omega\times \mathcal{I},\label{momentum}\\
 \partial_t\rho +   \nabla \cdot (\rho \bu) &= 0 &\quad \text{in} &\quad \Omega\times \mathcal{I},\label{mass}\\
 \partial_t\rho + \bu \cdot \nabla \rho &= 0 &\quad \text{in} &\quad \Omega\times \mathcal{I}, \label{incompressiblity}\\
  \mathbf{u}(\mathbf{x},0) &= \mathbf{u}_0(\mathbf{x})  & \quad \text{in} &\quad \Omega, \label{IC velocity}\\
 \rho(\mathbf{x},0) &= \rho_0(\mathbf{x})  & \quad \text{in} &\quad \Omega, \label{IC density} 
\end{alignat}
\end{subequations}
for the fluid velocity $\mathbf{u}:\Omega \rightarrow \mathbb{R}^d$, the pressure $p:\Omega \rightarrow \mathbb{R}$ and the density $\rho:\Omega \rightarrow \mathbb{R}$.
The problem is augmented with appropriate boundary conditions.
The equations (\ref{momentum})-(\ref{IC velocity}) describe the balance of momentum, the conservation of mass, the incompressiblity constraint and the initial conditions, respectively.
We denote with $\mathbf{x}\in \Omega$ the spatial parameter and with $t \in \mathcal{I}=(0,T)$ the time with end time $T>0$. 
The dynamic viscosity $\mu: \Omega \rightarrow \mathbb{R}^+$ depends on the density, i.e. $\mu = \mu(\rho)$. 
Furthermore, the body force is $\mathbf{g}: \Omega \times \mathcal{I} \rightarrow \mathbb{R}^d$ (this is often the gravitational force), the initial velocity 
is $\mathbf{u}_0: \Omega \rightarrow \mathbb{R}^d$ and the initial density is $\rho_0: \Omega \rightarrow \mathbb{R}$.
We assume a zero-average pressure for all $t \in \mathcal{I}$. 
The various derivative operators are the temporal one $\pd_t$ and the symmetric gradient $\nabla^s=\tfrac{1}{2}\left(\nabla + \nabla^T\right)$. The normal velocity is $u_n := \bu \cdot \mathbf{n}$.

In this paper we employ the level-set method to describe the two-fluid interface.
Hence, we define the scalar function $\phi=\phi(\mathbf{x},t)$ to distinguish the time-dependent subdomains of the two fluids, $\Omega^0_t$ and $\Omega^1_t$ respectively, via an interface $\Gamma_t$:
\begin{subequations}\label{domain decomposition}
\begin{alignat}{3}
 \Omega^0_t:= \left\{ \mathbf{x}\in \Omega ~|~ \phi(\mathbf{x},t)>0 \right\},\label{water}\\
 \Omega^1_t:= \left\{ \mathbf{x}\in \Omega ~|~ \phi(\mathbf{x},t)<0 \right\},\label{air}\\
 \Gamma_t:= \left\{ \mathbf{x}\in \Omega ~|~ \phi(\mathbf{x},t)=0 \right\}.\label{interface}
\end{alignat}
\end{subequations}
The fluid properties of the two subdomains are determined by the indicator function:
\begin{subequations}\label{fluid prop}
\begin{align}
\rho =& \rho_0 (1-H) + \rho_1 H,  \label{eq:rho_def} \\
\mu =& \mu_0 (1-H) + \mu_1 H.  \label{eq:mu_def} 
\end{align} 
\end{subequations}
The function $H=H(\phi)$ indicates the subdomain. The constant densities of two fluids are $\rho_0$ and $\rho_1$, and $\mu_0$ and $\mu_1$ are the constant dynamic viscosities of fluids.
For convenience we introduce the notation $\Delta \rho = \rho_1-\rho_0$ for the jump in density and $\Delta \mu = \mu_1-\mu_0$ for the jump in dynamic viscosity.
Note that a sharp interface requires the indicator function $H$ to be discontinuous.
A discrete implementation often uses a continuous indicator function.
This function gets arbitrarily close to the discontinuous indicator function in case of a vanishing transition domain.
In the following we assume a smooth $H$ to simplify the exposition.
We describe discrete implementation of $H=H(\phi)$ in subsection \ref{sec:LSER}.

By combining (\ref{mass})-(\ref{incompressiblity}) with (\ref{fluid prop}) the governing equations take the form:
\begin{subequations}\label{eq:strong_cons_ns}
\begin{alignat}{2}
 \partial_t (\rho \bu) + \nabla \cdot (\rho \bu\otimes \bu) 
 + \nabla  p -  \nabla \cdot 2\mu \nabla^s \bu &=  \rho \mathbf{g}, \label{momentum 2}\\
 \nabla \cdot \bu &= 0, \label{mass 2}\\
 \partial_t H + \bu \cdot \nabla H &= 0, \label{incompressiblity 2}
\end{alignat}
\end{subequations}
where we omit the initial condition and the specification of domains for convenience.

We assume no-penetration  boundary conditions, i.e. $\mathbf{n} \cdot \bu = 0$ on  $\Gamma$,
 and free-slip boundary conditions. 
Defining the appropriate velocity 
$\mathcal{U} =\{  \mathbf{u} \in  [H^1(\Omega) ]^d ;  \mathbf{u}\cdot\mathbf{n}=0 \} $
and pressure space
 $ \mathcal{P} =      \{p \in L^2(\Omega); \int p {\rm d}\Omega =0 \} $,
the  standard conservative weak formulation corresponding to the strong form
(\ref{eq:strong_cons_ns}) reads:\\



\textit{Find $\mathbf{u} \in \mathcal{U} ,~p \in \mathcal{P} ,~H \in H^1(\Omega)$ 
such that for all}
\textit{$~\mathbf{w} \in \mathcal{U},~q \in \mathcal{P} ,~\added[id=Rev3.2]{\psi \in L^2(\Omega)}$,}

\begin{subequations}\label{weak form 1}
\label{eq:weak}
\begin{align}
(\bw,  \partial_t(\rho \bu)) -(\nabla \bw, \rho \bu\otimes \bu) 
- (\nabla \cdot \bw, p) + (\nabla \bw, 2\mu \nabla^s \bu) 
 &= (\bw, \rho \mathbf{g}), \label{weak form 1 mom}\\
 ( q, \nabla \cdot \bu) &= 0, \label{weak form 1 cont} \\
(\psi, \partial_t H) + (  \psi,  \bu \cdot \nabla  H) &= 0, \label{weak form 1 mass}
\end{align}
\end{subequations}
where $(\cdot, \cdot)$ is the $L^2(\Omega)$ inner product on the interior. All boundary terms vanish due to the choice of boundary conditions.
Note that the weak formulation (\ref{weak form 1}) is equivalent to the strong form (\ref{eq:strong_cons_ns}) 
in the case of sufficiently \added[id=Rev2.1]{smooth} solutions.

This section continues with the discussion about the conservation of mass, momentum and energy associated with the weak formulation (\ref{weak form 1}).
This provides the blueprint on how to obtain conservation statements in the discrete setting, which is in its turn presented in section \ref{sec:disc_cons_form}.

\subsection{Conservation properties}

To derive conservation properties we select  appropriate weighting functions.
As we are still on the continuous level, we assume the selection of the appropriated weights  is allowed under reasonable restrictions.
This does not apply in the discrete setting, which is presented in section \ref{sec:disc_cons_form}.

\subsubsection{Mass} 

The conservation of mass follows when we select the weighting functions $\bw = 0$, $q=H\Delta \rho $ and $\psi=\Delta \rho$ in 
the weak statement (\ref{weak form 1}):
\begin{subequations}
 \begin{align}
 (H \Delta \rho , \nabla \cdot \bu) &= 0, \label{cons mass 1}\\
(\Delta \rho, \partial_t H) + ( \Delta \rho , \bu \cdot \nabla H) &= 0. \label{cons mass 2}
\end{align}
\end{subequations}
Next, we (i) apply Green's identity on (\ref{cons mass 2}), (ii) use that $\Delta \rho$ is constant and (iii) use that the domain does not change in time.
Combining the resulting equations delivers the conservation of global mass:
\begin{align}\label{eq:mass_cont}
\frac{{\rm d}}{{\rm d}t}\int_\Omega  \rho~  {\rm d}\Omega &= 0.
\end{align}

\subsubsection{Momentum} 
To show the conservation of linear momentum, the weak form (\ref{eq:weak}) needs to be augmented with an auxiliary flux \cite{HEML00, HuWe05}.
This is to be able to select the appropriate weighting function on the boundary $\Gamma$.
We note that this approach can be understood as a Lagrange multiplier construction \cite{EiAk17i, EiAk17ii}.

The problem takes the form:\\


\textit{Find $\mathbf{u} \in [H^1(\Omega)]^d,~p \in \mathcal{P},~H \in  H^1(\Omega)$ 
and $\lambda \in H^{-\frac{1}{2}}(\Gamma)$ such that for all}
~$\mathbf{w} \in [H^1(\Omega)]^d, ~q \in \mathcal{P},~\added[id=Rev3.2]{\psi \in  L_2(\Omega)}$,
and $\eta \in H^{-\frac{1}{2}}(\Gamma)$

\begin{subequations}\label{weak form 1 aux flux}
\begin{align}
(\bw,  \partial_t(\rho \bu)) -(\nabla \bw, \rho \bu\otimes \bu) 
- (\nabla \cdot \bw, p) + (\nabla \bw, 2\mu \nabla^s \bu) 
&= (\bw, \rho \mathbf{g}) + \langle \B{n}\cdot \bw,\lambda \rangle, \label{weak form 1 aux flux mom}\\
 ( q, \nabla \cdot \bu) &= 0, \label{weak form 1 aux flux cont}\\
(\psi, \partial_t H) + (  \psi,  \bu \cdot \nabla  H) &= 0 \label{weak form 1 aux flux mass} \\
\langle \eta,  \B{n}\cdot \bu \rangle &= 0 .\label{weak form 1  no pen} 
\end{align}
\end{subequations}
where $\langle \cdot , \cdot \rangle$ is a duality pairing on the boundary \footnote{\added[id=Rev3.2]{Note that $\mathbf{n}\cdot \bu \in H^{\frac{1}{2}}(\Gamma)$,  see \cite{Temam}.}}.
This formulation allows the choice of the weighting functions $\bw = \mathbf{e}_i$,  $q=0$ and $\psi=0$.
Substitution into (\ref{weak form 1 aux flux}) gives:
\begin{align}
(\mathbf{e}_i,  \partial_t(\rho \bu)) = (\mathbf{e}_i, \rho \mathbf{g}) + \langle  \B{n}\cdot  \mathbf{e}_i,\lambda\rangle,
\end{align}
where $\mathbf{e}_i$ are the Cartesian unit vectors. 
Equivalence of (\ref{weak form 1 aux flux}) with the strong form (\ref{eq:strong_cons_ns}) provides via Green's identity the expression of the auxiliary flux \cite{EiAk17i, EiAk17ii}:
\begin{align}
 \lambda = - p  + 2 \mu \mathbf{n} \cdot (\nabla^s \mathbf{u} ~\mathbf{n}).
\end{align}
The auxiliary flux is the generalized constraint force which enforces the no-penetration constraint.
Given that $\mathbf{e}_i$ is an arbitrary unit vector and assuming a time-independent domain, we arrive at:
\begin{align}\label{eq:mom_cont}
\frac{{\rm d}}{{\rm d}t}\int_\Omega  \rho \bu  \, {\rm d}\Omega &= \int_\Omega \rho \mathbf{g} {\rm d} \,\Omega - \int_\Gamma p \mathbf{n}  \, {\rm d}\Gamma  + \int_\Gamma  2\mu \, \mathbf{n}  \cdot  \nabla^s \bu \, {\rm d}\Gamma,
\end{align}
which is conservation of  momentum.

\subsubsection{Energy} 
The conservation of energy follows when selecting $\bw = \bu$ and  $q=p$ in the formulation (\ref{weak form 1}). The choice of $\psi$ is postponed. 
Assuming no-slip boundary conditions leads to
\begin{align} \label{eq:energy} 
(\bu,  \partial_t (\rho \bu)) -(\nabla \bu, \rho \bu\otimes \bu) 
+2 \| \mu^{1/2} \nabla^s \bu\|^2  &= (\bu, \rho \mathbf{g}).
\end{align}
We proceed by analyzing the kinetic and potential energy behavior separately. The evolution of the kinetic energy is linked to the first two terms on the left-hand side of (\ref{eq:energy}), while the evolution of the potential energy is associated with the right-hand side of (\ref{eq:energy}). Both require a different weighting function $\psi$.

\subsubsection*{Kinetic energy} 
\added[id=Rev2.2]{
The evolution of the kinetic energy is governed by the acceleration and
convective term in (\ref{eq:energy}).  To show this both terms are rewritten.
The acceleration term can be rewritten as
\begin{align} \label{eq:acc_term}
(\bu,  \partial_t (\rho \bu)) =&(\bu,  \bu \partial_t \rho ) +(\bu,  \rho \partial_t  \bu) \nonumber \\
=&\onehalf (\bu\cdot \bu,  \partial_t \rho ) + \frac{{\rm d}}{{\rm d}t} E_{kin}.
\end{align}
The kinetic energy is defined as
\begin{align}\label{eq:kin_eng}
E_{kin} := \int_\Omega \tfrac{1}{2}\rho \bu\cdot \bu \,{\rm d}\Omega.
\end{align}
We assume that the domain does not change in time.
The convective term can be rewritten as
\begin{align} \label{eq:conv_term}
  - \left(\nabla \mathbf{u}, \rho \mathbf{u}\otimes \mathbf{u} \right) =& \frac{1}{2} \left( \mathbf{u},\nabla \left( \rho \mathbf{u}\otimes \mathbf{u}\right) \right) - \frac{1}{2} \left(\nabla \mathbf{u}, \rho \mathbf{u}\otimes \mathbf{u} \right)\nn\\
								       =& \frac{1}{2} \left( \mathbf{u},\mathbf{u} (\mathbf{u} \cdot \nabla \rho ) \right) +\frac{1}{2} \left( \mathbf{u},\rho \mathbf{u}\cdot \nabla \mathbf{u}  \right) +  \frac{1}{2} \left( \mathbf{u},\rho \mathbf{u} \nabla \cdot  \mathbf{u} \right)- \frac{1}{2} \left(\nabla \mathbf{u}, \rho \mathbf{u}\otimes \mathbf{u} \right)\nn\\
								       =& \frac{1}{2} \left( \mathbf{u},\mathbf{u} (\mathbf{u} \cdot \nabla \rho ) \right) +  \frac{1}{2} \left( \mathbf{u},\rho \mathbf{u} \nabla \cdot  \mathbf{u} \right).
\end{align}
At this point it is important to emphasize that the convective contribution does not vanish, even for divergence-free velocity fields.
This is in contrast to the single fluid case where it does vanish, see e.g. \cite{palha2017mass, lee2017discrete, EiAk17ii}.
The difference lies in the presence of the density.
In general the density varies between the two fluids and thus its gradient is not identically zero.

Note that the volume term with a density gradient reduces to an interface term with a density jump in case of a vanishing 
smoothing distance. }

A straightforward combination of eq (\ref{eq:acc_term}) and (\ref{eq:conv_term}) reveals:
\begin{align}\label{intermediate result zoveel}
(\bu,  \partial_t (\rho \bu))  - (\nabla \bu, \rho \bu\otimes \bu) = & \frac{{\rm d}}{{\rm d}t} E_{kin} +  (\partial_t \rho + \bu \cdot \nabla \rho, \tfrac{1}{2} \bu\cdot \bu)+  \frac{1}{2} \left( \mathbf{u},\rho \mathbf{u} \nabla \cdot  \mathbf{u} \right),
\end{align}

We select $\psi = \tfrac{1}{2} \Delta \rho  \bu \cdot \bu$ in (\ref{weak form 1 mass}) and modify $q = p + \tfrac{1}{2}\rho \bu\cdot\bu$, recall  (\ref{eq:rho_def}) and add the result to (\ref{intermediate result zoveel}). This causes the second and third term on the right-hand side to vanish and we are left with:
\begin{align}\label{eq:kin_energy}
\frac{{\rm d}}{{\rm d}t} E_{kin} & = (\bu, \partial_t (\rho \bu))  - (\nabla \bu, \rho \bu\otimes \bu),
\end{align}
which is the asserted relation alluded earlier.

\subsubsection*{Potential energy} 
The evolution of the potential energy is governed by the body force term on the left-hand side of  (\ref{eq:energy}). 
Assuming no-slip boundary conditions this term can be written as:
\begin{align}\label{cont pot}
(\bu, \rho \mathbf{g}) 
&= (\bu,  (\nabla \bx) \rho \mathbf{g})  \nn \\
&= 
 - (\nabla \cdot (\rho \bu) ,    \bx \cdot \mathbf{g})   \nn \\
&=  
(\partial_t \rho,  \bx \cdot \mathbf{g})   - (\bx \cdot \mathbf{g} ,  \partial_t \rho + \nabla \cdot (\rho \bu) ).
\end{align}
The last term on the right-hand side vanishes if we select $\psi = \Delta \rho\,  \bx \cdot \mathbf{g} $ and recall  (\ref{eq:rho_def}).
Using the definition for the potential energy
\begin{align} \label{eq:pot_eng}
E_{pot} := - \int_\Omega \rho \bx \cdot \mathbf{g}   \,{\rm d}\Omega, 
\end{align}
and assuming that the domain does not change in time, we arrive at:
\begin{align}\label{eq:pot_energy}
\frac{{\rm d}}{{\rm d}t} E_{pot} &=  -   (\bu, \rho \mathbf{g}),
\end{align}
which proves the claim.

\subsubsection*{Total energy}

We refocus our attention on the overall energy evolution relation as stated in (\ref{eq:energy}).
Using the derived equations relating to 
kinetic energy (\ref{eq:kin_energy})
and
potential energy (\ref{eq:pot_energy}), the overall energy evolution takes the form:
\begin{align}\label{eq:eng_cont}
\frac{{\rm d}}{{\rm d}t} E_{kin} + 2 \| \mu^{1/2} \nabla^s \bu\|^2 
&=  -\frac{{\rm d}}{{\rm d}t} E_{pot},
\end{align}
which is the conservation of energy.
This clearly states that in the inviscid case, kinetic and potential energy are exchanged, and the total energy is therefore conserved.\\


\section{Standard conservative discretization}
\label{sec:disc_cons_form}
This section introduces a discrete formulation that closely resembles the continuous formulation of the previous section. 
We examine the conservation properties of this formulation. 
The Crank-Nicolson method is employed for the temporal discretization, 
while a combination of Galerkin, SUPG and level--sets is employed for the spatial discretization.
We start this section by providing a brief review of isogeometric analysis which serves as an important concept of the developed methodology.
Next, we discuss the level-set method and subsequently the standard discretization with its conservation properties.

\subsection{Isogeometric Analysis}

\label{sec:space}
To discretize the governing fluid flow equations we use the isogeometric analysis (IGA) technology developed by Hughes and coworkers in \cite{HuCoBa04}.
IGA integrates the historically different fields of Computer-Aided Design (CAD) and Computer-Aided Engineering (CAE).
It unifies the representation of the geometry of the CAD design and CAE analysis. 
To this purpose NURBS (Non-Uniform Rational B-Spline) basis functions are employed.
For the formal definition of these shape functions we refer to \cite{CoHuBa09} and \cite{PiegTil97}.
The parametrization of the solution is the same as that of the underlying geometry, which is known as the isoparameteric concept.
Isogeometric analysis shares many features in common with finite elements, such as the underlying variational
framework, the isoparametric concept, locally supported basis functions, and possibilities for $h$- and $p$-refinement.
However, isogeometric analysis offers richer possibilities for geometry modeling and solution representation, compared to the finite element methods.
For example, the NURBS surfaces in IGA match \textit{exactly} the CAD geometry and IGA offers an improved refinement strategy known as $k$-refinement, 
which is not possible in FEM.
Due to this $k$-refinement the global degree-of-freedom count is smaller in comparison with standard finite elements.
This leads to efficient higher-order discretizations, both in theory \cite {EvBaBaHu08, BBCHS06} and in applications, see for instance \cite {BACHH07}. 

In this paper the geometric features are not of direct importance.
Here we take advantage of the higher-order and higher-continuity properties of the NURBS.
The ability to control the interelement continuity, besides the element order itself, gives the 
NURBS the flexibility to construct combinations of velocity and pressure discretization inconceivable before.
Well-known finite-element families, such as Taylor-Hood elements, are defined for different orders. 
These families can be extended by adding interelement continuity as a new parameter.
In \cite{Evans13steadyNS, Evans13unsteadyNS, buffa2011isogeometric, buffa2011isogeometric3D} this 
additional parameter is used to construct stable  velocity and pressure pairs that allow 
\textit{pointwise}  divergence-free velocity fields. 
In the subsection \ref{sec:cons_form} we show that the discretization indeed employs pointwise solenoidal velocities.


\subsection{Level--set approach with explicit redistancing}\label{sec:LSER}
Here we include a brief description of the level-set method with explicit redistancing proposed in \cite{AKKERMAN2017}.

To describe the fluid interface the indicator function $H=H(\phi)$ could be represented by a simple Heaviside function:
\begin{align}\label{eq:heavi}
H(\phi) = \left \{
\begin{array}{lcr}
0& \text{if} &\phi < 0,\\
\onehalf & \text{if} &\phi =0,\\
1& \text{if} &\phi > 0,
\end{array}
\right . 
\end{align}
where the positive and negative part each represent one of the two fluids and the zero-level is the interface.
This definition leads to problems when directly employed in a numerical method. 
Therefore, the sharp Heaviside function is often smoothed.
There are several option to perform the smoothing. 
A popular choice, that we also adopt here, is to take as smoothed Heaviside function:
 \begin{align}\label{eq:heavi_smooth}
\hat{H}(\phi_\alpha) = \left \{
\begin{array}{lcr}
0& \text{if} &\phi_\alpha \leq -1,\\
\onehalf(1 + \sin(\frac{\pi}{2} \phi_\alpha)) & \text{if} & |\phi_\alpha| < 1,\\
1& \text{if} &\phi_\alpha \geq 1,
\end{array}
\right . 
\end{align}
where $\phi_\alpha$ is a scaled level-set.

Scaling and redistancing techniques of the level-set are necessary to properly control the smoothing region around the interface.
Traditionally, the redistancing step is done by solving the Eikonal-equation \cite{Sethian99,Sethian_01}.
This is a nonlinear problem which makes it hard to include it in a monolithic solver.
The demanded computational effort and lack of robustness are often the main concerns.
Another issue involves the required large number of iterations to arrive at a redistanced level-set.
An additional complication is the trade-off between on the one hand the actual redistancing around the interface and on the other maintaining the interface location.
These requirements contradict and therefore a compromise is demanded.

The scaling is often done by only taking the local mesh size into account.
This can cause problems in highly graded meshes.
Basing the scaling on the average mesh size between the current point and the closest point to the interface, circumvents this issue. 
Note that this would require the scaling to be based on an integral quantity.

These redistancing and scaling issues have been tackled in \cite{AKKERMAN2017}. 
The crucial step is to introduce a scaling parameter $\alpha$ that relates the convected $\phi$ with the 
redistanced and rescaled $\phi_\alpha$ via $\phi_\alpha=\phi/\alpha$.
This directly solves the paradox of redistancing and maintaining the interface.
Independent of the scaling, the zero-level-set of $\phi_\alpha$ and $\phi$ are identical.
It turns out that both redistancing and scaling can be achieved by solving a simple projection for $\alpha$:\\

\textit{Find $\alpha \in H^1(\Omega)$ such that for all $\eta \in H^1(\Omega)$,}
\begin{align} \label{eq:alpha1}
\left (  \eta,   \| \nabla_\xi \phi \|  \right ) = \left (  \eta,   \alpha \right )  + \epsilon \left (  \nabla_\xi \eta,    \nabla_\xi \alpha \right ).
\end{align}
Here  $\epsilon$ is a given smoothing parameter and $\nabla_\xi$ is the gradient with respect to the reference coordinate $\xi$.

\subsection{Discrete weak formulation}\label{sec:cons_form}
The weak formulation employs for a large part the standard discretization of the continuous conservative form (\ref{weak form 1}).
The discretization of the momentum and continuity equations indeed uses the staightforward Galerkin method.
Hence, the methodology is only applicable for the computation of low Reynolds number flows; the potential convective instability that can occur because of the Galerkin discretization is thus circumvented.
The use of stabilization techniques to deal with high Reynolds flows is beyond the scope of this paper.
However, the spatial discretization of the interface evolution equation does require stabilization.
This equation is a pure convection problem and therefore a standard discretization is prone to wiggles.
Accordingly, we adopt the well-known SUPG method \cite{BroHug82} as a stable discretization of the interface convection equation.
Furthermore, since the level-set methodology is employed, the level-set $\phi$ replaces the indicator function $H$. 

\subsection*{Remark}
The Galerkin method requires the use of compatible velocity and pressure discretizations.
The isogeometric analysis concept, which employs NURBS basis functions, is employed for this purpose.
The NURBS-spaces guarantee exact divergence-free velocity fields.
Subsection \ref{sec:space} provides some background on this.\\

We employ the Crank-Nicolson method for the temporal integration.
This method is an unconditionally stable second-order integrator.
The motivation for this choice emerges from an energy perspective.
In a mono-fluid setting, i.e. a constant density, this method is an energy conservative 
time-integrator.
In fact, within a generalized-$\alpha$ framework, it is the only second-order method that can be linked to proper energy decay.
This is in detail described in \cite{EiAk17i}.

Employing the described ingredients we arrive at the discrete weak formulation:\\
\textit{Find $\mathbf{u}^{n+1} \in \mathcal{U} ,~p^{n+1} \in \mathcal{P} ,~\phi^{n+1} \in H^1(\Omega)$ 
such that for all}
\textit{$~\mathbf{w} \in \mathcal{U},~q \in \mathcal{P} ,~\psi \in H^1(\Omega)$,}

\begin{subequations}\label{eq:weak form}
\label{eq:cons_form}
\begin{align}
\left (\bw,  \frac{\rho^{n+1} \bu^{n+1} - \rho^{n} \bu^{n}}{\Delta t} \right) -(\nabla \bw, \rho^{n+1/2} \bu^{n+1/2} \otimes \bu^{n+1/2} ) \nn \\
- (\nabla \cdot \bw, p^{n+1} ) + (\nabla \bw, 2\mu \nabla^s \bu^{n+1/2} ) &= (\bw, \rho^{n+1/2}  \B{g}), \\
 ( q, \nabla \cdot \bu^{n+1/2} ) &= 0, \label{incompressibility discrete}\\
\left (\psi, \frac{\phi^{n+1} - \phi^{n} }{\Delta t} + \bu^{n+1/2}  \cdot \nabla \phi^{n+1/2} \right)  \nn\\
+\left (  \tau \bu^{n+1/2}  \cdot \nabla \psi,  \frac{\phi^{n+1} - \phi^{n} }{\Delta t} + \bu^{n+1/2}  \cdot \nabla \phi^{n+1/2} \right ) 
&= 0,\label{discrete interface}
\end{align}
\end{subequations}
\added[id=Rev3.2]{where the test space $\psi$  is reduced to $H^1(\Omega)$ to accommodate the SUPG terms.}
The stabilization parameter is defined as in \cite{ABKF11}: 
\begin{align}\label{eq:tau_def}
\tau = \left ( \frac{4}{\Delta t^2} + \bu^{n+1/2} \cdot \mathbf{G}  \bu^{n+1/2} \right )^{-1/2},
\end{align}
and where $\mathbf{G}$ is the second-rank metric tensor:
\begin{align}
\mathbf{G}   = \frac{\partial \boldsymbol{\xi}} {\partial \bx }^T \frac{\partial \boldsymbol{\xi}} {\partial \bx }.
\end{align}
Here $\mathbf{x}$ and $\boldsymbol{\xi}$ are the spatial coordinates in the physical and parameter domain respectively.
The superscripts $n, n+1/2, n+1$ refer to the current, intermediate and next time-level, respectively, 
and $\Delta t$ denotes the time step. The fluid properties are given by
\begin{subequations}
 \label{eq:mu_rho_def2}
\begin{align}
\rho =& \rho_0 \left (1-\hat{H}\left (\frac{\phi}{\alpha} \right ) \right )+ \rho_1 \hat{H}\left (\frac{\phi}{\alpha} \right ),   \label{eq:rho_def2} \\
\mu =& \mu_0 \left (1-\hat{H}\left (\frac{\phi}{\alpha} \right ) \right ) + \mu_1 \hat{H}\left (\frac{\phi}{\alpha} \right ),   \label{eq:mu_def2}
\end{align} 
\end{subequations}
where $\alpha $ is computed via (\ref{eq:alpha1}).

Note that the choice of the time-level in (\ref{incompressibility discrete}) is not standard. It is required to arrive at the desired energy behavior, see section \ref{sec:cons_eng}.

\subsubsection{Divergence-free solutions}\label{sec:div_free}

As mentioned in subsection \ref{sec:space} and detailed in  \cite{Evans13steadyNS, Evans13unsteadyNS, buffa2011isogeometric, buffa2011isogeometric3D} 
the ability to tune the interelement continuity of the NURBS functions allows the construction of favorable velocity and pressure spaces.
In particular they can be chosen such that the divergence of the velocity is a member of the pressure space.
On general non-aligned or curved meshes,  this requires the use of the Piola transformation \cite{arnold2005quadrilateral, rognes2009efficient}.

This allows us to prove that the velocities are point-wise divergence-free
by selecting the weights $\bw = 0$,  $q=\nabla \cdot \bu^{n+1/2} $ and $\psi=0$ in (\ref{eq:weak}). It is important to realize that this choice is only possible due to the delicate choice of the discretization spaces. We arrive at
\begin{align}
  \| \nabla \cdot \bu^{n+1/2}\|^2 = 0,
\end{align}
this directly leads to,
\begin{align}
\nabla \cdot \bu^{n+1/2} = 0, \quad \text{for all}~x \in \Omega.
\end{align}
This means that the solution is point-wise divergence-free at the midpoint in time.
Using the definition for the midpoint velocity vector we get:
\begin{align}
\nabla \cdot \bu^{n+1} = - \nabla \cdot \bu^{n}. 
\end{align}
This means that the divergence error from one time step is directly mirrored into the next time step. 
Hence, the discretization provides pointwise divergence-free solutions for a solenoidal initial condition.

Alternatively, the velocity in the continuity equation can be based on time-step $n+1$. 
This would lead to divergence-free solutions independent of the initial condition.
However, in this case the discretized version of the energy statement given in (\ref{eq:energy}) would be augmented with the term
\begin{align}
\onehalf (\nabla \cdot \bu^{n+1}  - \nabla \cdot \bu^{n}, p^{n+1}),
\end{align}
due to the mismatch in time-levels of the velocities.
When employing point-wise divergence-free solutions this term vanishes (except for possibly the first time-step).


\subsection{Conservation properties}
\subsubsection{Mass} \label{sec:cons_mass}
The straightforward choice of weighting functions to arrive at a statement about the conservation of mass 
would be to take the discrete counterpart of the continuous case, i.e. $\bw = 0$,  $q=\Delta \rho \hat{H} $ and $\psi=\Delta \rho \partial \hat {H}/\partial \phi$. 
However, this is not a valid choice: both the functions $\hat{H}$ and $\partial \hat{H}/\partial \phi$ are not in the weighting function spaces.
Mass conservation can therefore not be guaranteed {\it a priori}.
To equip the discrete formulation with the mass conservation property, it should be explicitly enforced.
This means that the constraint
\begin{align}\label{eq:disc_mass_constr}
h_{1}(\phi^{n+1})  := \int_\Omega ( \rho^{n+1} - \rho^n) ~  {\rm d}\Omega &= 0
\end{align}
should be fulfilled to guarantee global mass conservation.

\subsubsection{Momentum} \label{sec:cons_mom}
In contrast to the conservation of mass, for the conservation of momentum the same procedure as in the continuous case can be directly employed. 
This results in the global conservation statement:
\begin{align}
\Delta t^{-1}\int_\Omega ( \rho^{n+1} \bu^{n+1}  -\rho^{n}  \bu^{n})   \, {\rm d}\Omega &= 
\int_\Omega \rho^{n+1/2} \B{g} {\rm d} \,\Omega 
- \int_\Gamma p^{n+1}  \B{n}  \, {\rm d}\Gamma  
+ \int_\Gamma  2\mu \, \B{n}  \cdot  \nabla^s \bu^{n+1/2}  \, {\rm d}\Gamma,
\end{align}
which is the straightforward discrete counterpart of (\ref{eq:eng_cont}).

\subsubsection{Energy} \label{sec:cons_eng}
The statements about the conservation of the discrete energies follow when the discrete 
counterparts of the continuous weights are chosen: $\bw = \bu^{n+1/2}$ and  $q=p^{n+1}$. 
The specific choice for $\psi$ is again postponed.
This  leads to the discrete equivalent of (\ref{eq:energy}) which states:
\begin{align}
\left ( \bu^{n+1/2},  \frac{\rho^{n+1} \bu^{n+1} - \rho^{n} \bu^{n}}{\Delta t} \right)
 -\left(\nabla  \bu^{n+1/2}, \rho^{n+1/2} \bu^{n+1/2} \otimes \bu^{n+1/2} \right) \nn \\
+~ 2 \| \mu^{1/2} \nabla^s \bu^{n+1/2} \|^2 &= \left( \bu^{n+1/2}, \rho^{n+1/2}  \B{g}\right).
\end{align}
We proceed in a similar way as in the continuous case.

\subsubsection*{Kinetic energy}\label{sec:cons_kin}
The discrete acceleration and convective terms can be expressed as:
\begin{align}
\left ( \bu^{n+1/2},  \frac{\rho^{n+1} \bu^{n+1} - \rho^{n} \bu^{n}}{\Delta t} \right) 
-(\nabla  \bu^{n+1/2}, \rho^{n+1/2} \bu^{n+1/2} \otimes \bu^{n+1/2} ) 
= \nn \\
\frac{E_{kin}^{n+1} - E_{kin}^{n}}{\Delta t}+\left (\frac{\rho^{n+1} - \rho^{n}}{\Delta t} ,\tfrac{1}{2}\bu^{n}\cdot \bu^{n+1} \right )+\left (\nabla \cdot \left(\rho^{n+1/2} \bu^{n+1/2}\right),\tfrac{1}{2}\bu^{n+1/2}\cdot \bu^{n+1/2} \right ),
\end{align}
where the kinetic energies are defined by taking (\ref{eq:kin_eng}) at the corresponding time-level. This is a discrete version of (\ref{intermediate result zoveel}).
 In the continuous case the second term on the right-hand side cancels by choosing an appropriate weight in the 
 interface evolution equation (\ref{discrete interface}). 
 Similarly to the mass conservation, the density needs to be related to the level--set by select a weight proportional to 
$H$ or $\frac{\partial H}{\partial \phi}$ which is not allowed.
Additionally, (i) the time levels of the velocity fields  in the time derivative and convection term do not match 
and (ii) the unwanted terms caused by the SUPG stabilization pollute the relation. To ensure a link between the acceleration and convective term with the kinetic energy, the constraint
\begin{align}\label{eq:disc_kin_eng_constr}
h_{2}(\phi^{n+1}) := \left (\frac{\rho^{n+1} - \rho^{n}}{\Delta t},\tfrac{1}{2}\bu^{n}\cdot \bu^{n+1} \right ) 
- (\rho^{n+1/2} \bu^{n+1/2}, \bu^{n+1/2} \cdot \nabla \bu^{n+1/2})  =0
\end{align}
needs to be explicitly enforced. Note that Green's identity has been applied in the last term.
The divergence-free velocities provided by the isogeometric analysis framework have positive benefits for the constraints.
Employing this property we can formulate the constraint as
\begin{align}\label{eq:disc_kin_eng_constr 2}
\left (\frac{\rho^{n+1} - \rho^{n}}{\Delta t},\tfrac{1}{2}\bu^{n}\cdot \bu^{n+1} \right ) 
+ (\bu^{n+1/2}\cdot \nabla \rho^{n+1/2},\tfrac{1}{2}\bu^{n+1/2}\cdot \bu^{n+1/2})  =0.
\end{align}
Here we recognize a convective interface contribution in the second term.
Let us now consider the mono-fluid case, i.e. $\rho_0 = \rho_1$ and $\mu_0 = \mu_1$. 
In this setting the convection equation would be superfluous.
The constraint takes the form: 
\begin{subequations}
\begin{align}
 (\tfrac{1}{2}  \bu^{n+1/2} \cdot \bu^{n+1/2} , \nabla \cdot \bu^{n+1/2}  )=0,
\end{align}
\end{subequations}
which is obviously fulfilled when dealing with solenoidal velocity fields.
In other words, in the case of a constant density, the use of a velocity-pressure pair 
that results in divergence-free velocities is essential for the 
formulation stated in (\ref{eq:ce_form}) to yield an energy conservation statement.

\subsubsection*{Potential energy}\label{sec:cons_pot}
The discrete counterpart of (\ref{cont pot}) is:
\begin{align}\label{discrete pot}
( \bu^{n+1/2}, \rho^{n+1/2}  \B{g}) 
=
-  \frac{E^{n+1}_{pot}-E^{n}_{pot}}{\Delta t}- \left ( \frac{\rho^{n+1} - \rho^n}{\Delta t} +  \nabla \cdot (\rho^{n+1/2} \bu^{n+1/2}),   \bx \cdot  \B{g} \right),
\end{align}
where the potential energies are defined by taking (\ref{eq:pot_eng}) at the corresponding time level.
For reasons similar as before, a weight $\psi$ can not be chosen to ensure that the last term of (\ref{discrete pot}) vanishes.
Hence, we need to enforce the constraint
\begin{align}\label{eq:disc_pot_eng_constr}
h_{3}(\phi^{n+1}) :=  \left ( \frac{\rho^{n+1} - \rho^n}{\Delta t},   \bx \cdot  \B{g} \right)-  ( \rho^{n+1/2} , \bu^{n+1/2} \cdot  \B{g}) = 0
\end{align}
to guarantee a direct link between the body force and the potential energy. Again Green's identity has been applied in the last term.

Similar to the kinetic energy case, the divergence-free velocities yield favorable properties here. It reduces the constraint here to:
\begin{align}
 \left ( \frac{\rho^{n+1} - \rho^n}{\Delta t} +  \bu^{n+1/2} \cdot \nabla \rho^{n+1/2},   \bx \cdot  \B{g}\right)=0,
\end{align}
in which we again recognize a convective interface contribution in the second term.

\subsubsection*{Total energy}
Combining the previous results leads to the energy statement:
\begin{align}
\frac{E^{n+1}_{kin}-E^{n}_{kin}}{\Delta t} + h_{2}(\phi^{n+1})
+ 2 \| \mu^{1/2} \nabla^s \bu^{n+1/2} \|^2 &=  - \frac{E^{n+1}_{pot}-E^{n}_{pot}}{\Delta t} -h_{3}(\phi^{n+1}),
\end{align}
which is the discrete counterpart of (\ref{eq:eng_cont}) with two additional terms, namely $h_{2}(\phi^{n+1})$ and $h_{3}(\phi^{n+1})$. 
The sign and magnitude of these terms is undetermined and therefore artificial energy growth can not be precluded.
In order to guarantee correct energy behavior the additional terms should vanish.\\

This paper proceeds in section \ref{sec:EC_form} by enforcing the constraints via a Lagrange multiplier approach.
This leads to a numerical method with correct energy behavior for solving two-fluid flow.

\section{The discrete energy-corrected formulation}
\label{sec:EC_form}

In this section we present the corrected version of the standard discretization (\ref{eq:weak form}) 
that satisfies the global conservation of mass and energy. 
Therefore we employ the Lagrange multiplier method to enforce the constraints obtained in section \ref{sec:disc_cons_form}.

First, we present a small sketch of this method in a general abstract setting, after which we apply 
this approach to the standard discretization (\ref{eq:weak form}). We close with the solution strategy of the discretized system.
\subsection{The Lagrange multiplier method in a general PDE setting} \label{sec:lm} 
Here we present a brief description of the Lagrange multiplier method in a general PDE setting. 
Let $\mathcal{V}$ be a suitable function space with the $L^2$-innerproduct $(\cdot, \cdot)$ and induced norm $\|\cdot\|$. 
Consider the constrainted problem for a linear operator $\mathscr{L}$ and the functionals $f,h \in L^2(\Omega)$:\\

\textit{Find $v \in \mathcal{V}$ such that,}
\begin{subequations}\label{constrained problem}
  \begin{alignat}{2}
     \mathscr{L}v=&f, \label{eq:lm_base}\\
h(v) =& 0,\label{eq:lm_constr}
     \end{alignat}
\indent \textit{where}
\begin{alignat}{2}
h(v) =& \int_\Omega \tilde{h}(v) {\rm d} \Omega.\label{glob constr}
  \end{alignat}
\end{subequations}
Note that (\ref{eq:lm_constr})-(\ref{glob constr}) represents the enforcement of a global constraint.
The standard variational formulation corresponding to (\ref{constrained problem}) is:\\

\textit {Find $v \in \mathcal{V}$  and $\lambda \in \mathbb{R}$ such that  for all $w \in \mathcal{V}$}
\begin{subequations}
\label{eq:lm_form}
\begin{align}
 \left(w ,\mathscr{L}v - f\right)
 + \lambda  \left ( w, \frac{ \partial  \tilde{h}}{\partial v} \right ) =&0, \label{eq:lm_form1}\\
h(v) =& 0.
\end{align}
\end{subequations}
Note that (\ref{eq:lm_form1}) can be converted into a strong form similar to (\ref{eq:lm_base}) augmented with a perturbation term that scales with the Lagrange multiplier $\lambda$. 
This term creates the freedom in (\ref{eq:lm_form1}) in order to satisfy the constraint (\ref{eq:lm_constr}).

Problem (\ref{eq:lm_form}) could either be solved directly, or via a procedure that circumvents the saddle point nature of the problem.
The latter approach splits the solution into two components of which the second one scales with $\lambda$:
\begin{align} \label{eq:sol_decomp}
v = v_f  + \lambda v_{\lambda}.
\end{align}
This decouples the problems into one purely linked to the PDE:\\

\textit {Find $v_f \in \WW$  such that  for all $w \in \WW$}
\begin{align}
 (w ,\mathscr{L}v_f - f) =&0, 
\end{align}
and one involving the perturbation:\\

\textit {Find $v_{\lambda} \in \WW$  such that  for all $w \in \WW$}
\begin{align}
 (w ,\mathscr{L}v_{\lambda}) 
 + \left ( w, \frac{ \partial  \tilde{h}}{\partial v} \right ) =&0 .
\end{align}
The Lagrange multiplier $\lambda$ follows from the constraint:
\begin{align}
 h(v_f  + \lambda v_{\lambda}) =  0. 
\end{align}
Note that the solution (\ref{eq:sol_decomp}) satisfies the weak form (\ref{eq:lm_form}) only if $\mathscr{L}$ is linear.

\subsection{The numerical formulation}
To enforce the constraints (\ref{eq:disc_mass_constr}), (\ref{eq:disc_kin_eng_constr}) and (\ref{eq:disc_pot_eng_constr}) in the weak formulation (\ref{eq:weak form}) 
we apply the methodology presented in the previous subsection to the interface evolution equation (\ref{discrete interface}).
In principle one could also choose to add the constraints to the momentum equation.
However, this would effect the energy behavior of the formulation, which is the primary quantity of interest.
Additionally, the analysis of the energy behavior in the continuous and discrete setting, presented in sections \ref{sec:continuous} and \ref{sec:disc_cons_form}, indicate that the correct evolution of the interface leads to the correct energy behavior. 
It is therefore natural to perturb the convection equations. 

Augmenting the formulation (\ref{eq:cons_form}) with the mass (\ref{eq:disc_mass_constr}), kinetic energy (\ref{eq:disc_kin_eng_constr}) and potential energy constraint (\ref{eq:disc_pot_eng_constr}) and perturbing the convection equation appropriately, we arrive at the energy-corrected formulation:\\

\textit{Find $\mathbf{u}^{n+1} \in \mathcal{U} ,~p^{n+1} \in \mathcal{P} ,~\phi^{n+1} \in H^1(\Omega)$ 
and $\lambda_i \in   \mathbb{R}^3$
such that for all}
\textit{$~\mathbf{w} \in \mathcal{U},~q \in \mathcal{P} ,~\psi \in H^1(\Omega)$,}


\begin{subequations}
\label{eq:ce_form}
\begin{align}
\left (\bw,  \frac{\rho^{n+1} \bu^{n+1} - \rho^{n} \bu^{n}}{\Delta t} \right) -(\nabla \bw, \rho^{n+1/2} \bu^{n+1/2} \otimes \bu^{n+1/2} ) \nn \\
- (\nabla \cdot \bw, p^{n+1/2} ) + (\nabla \bw, 2\mu \nabla^s \bu^{n+1/2} ) &= (\bw, \rho^{n+1/2}  \B{g}), \\
 ( q, \nabla \cdot \bu^{n+1/2} ) &= 0, \\ \nn\\
\left (\psi, \frac{\phi^{n+1} - \phi^{n} }{\Delta t}  + \bu^{n+1/2}  \cdot \nabla \phi^{n+1/2} \right)  \nn\\
+\left (  \tau \bu^{n+1/2}  \cdot \nabla \psi,  \frac{\phi^{n+1} - \phi^{n} }{\Delta t} + \bu^{n+1/2}  \cdot \nabla \phi^{n+1/2} \right )+\sum_{i=1,2,3} \lambda_{i} \delta h_{i} &=0, \\
 h_{i} &=0,   \quad {i=1,2,3},\label{eq:ce_form const} 
\end{align}
\textit{where}
\begin{align}
\delta h_{1} =& \left (\psi, \partial_\phi \rho\ \right ),  \label{eq:const_vars 1} \\
\delta h_{2} =& \Delta t^{-1}\left (\partial_\phi \rho\, \psi ,\bu^{n}\cdot \bu^{n+1} \right )
                          -  ( \bu^{n+1/2}  \cdot  \nabla \bu^{n+1/2}   ,  \bu^{n+1/2} \partial_\phi \rho\, \psi ),  \label{eq:const_vars 2}\\
\delta h_{3}  =& \Delta t^{-1} \left ( \partial_\phi \rho\,\psi ,   \bx \cdot  \B{g} \right) 
                         -  \onehalf(\partial_\phi \rho\, \psi, \bu^{n+1/2} \cdot  \B{g}). \label{eq:const_vars 3}
\end{align}
\end{subequations}
The constraints $h_{1},h_{2}$ and $h_{3}$ are given in (\ref{eq:disc_mass_constr}), (\ref{eq:disc_kin_eng_constr}) and (\ref{eq:disc_pot_eng_constr}), 
respectively, while $\delta h_{1},\delta h_{2}$ and $\delta h_{3}$ are their variations with respect to  the level-set function $\phi^{n+1}$.
The derivative of the density with respect to the level-set function can be computed as:
\begin{align}
\frac{\partial \rho}{\partial \phi} 
=  (\rho_1-\rho_0 )\frac{\partial \hat{H}}{\partial \phi_\alpha } \frac{\partial \phi_\alpha }{\partial \phi } 
=  (\rho_1-\rho_0 ) \frac{\partial \hat{H}}{\partial \phi_\alpha } \alpha  
\end{align}
where the derivative of the smooth Heaviside (\ref{eq:heavi_smooth}) is:
 \begin{align}\label{eq:dirac_smooth}
\frac{\partial \hat{H}}{\partial \phi_\alpha } = \left \{
\begin{array}{lcr}
0& \text{if} &\phi_\alpha < -1,\\
\frac{\pi}{4} \cos(\frac{\pi}{2} \phi_\alpha) & \text{if} & |\phi_\alpha|  \leq 1,\\
0& \text{if} &\phi_\alpha > 1.
\end{array}
\right . 
\end{align}
This is a smoothed version of the Dirac function. The stabilization parameter is defined in (\ref{eq:tau_def}) while the fluid parameters are determined via (\ref{eq:mu_rho_def2}) and (\ref{eq:alpha1}).

Due to the explicit enforcement of the constraints (\ref{eq:ce_form const}), the global mass and energy conservation are restored, as discussed in section \ref{sec:disc_cons_form}. 
The formulation (\ref{eq:ce_form}) obeys the following energy balance:
\begin{align}\label{eq:disc_eng_ce}
\frac{E^{n+1}_{kin}-E^{n}_{kin}}{\Delta t} 
+ 2 \| \mu^{1/2} \nabla^s \bu^{n+1/2} \|^2 &=  - \frac{E^{n+1}_{pot}-E^{n}_{pot}}{\Delta t},
\end{align}
which directly mirrors the energy balance of the continuos formulation (\ref{eq:eng_cont}).

\subsection{Solution strategy of the discrete system}
Here we describe our strategy to solve discrete system resulting from (\ref{eq:ce_form}).
\subsubsection{Matrix structure}
The formulation (\ref{eq:ce_form}) results in a slightly unusual structure of the problem due to the constraints.
A straightforward Newton linearization namely leads to a discrete system with the following block structure:
\begin {align}
\left [
\begin{array}{cccccc}
A      &G &B_1 & 0 &0 &0\\
G^T & 0& 0      & 0&0 &0\\
B_2&0&C&t_1 & t_2 & t_3\\
r_1^T&0&s_1^T&0&0 &0\\
r_2^T&0&s_2^T&0&0 &0\\
r_3^T&0&s_3^T&0&0 &0\\
\end{array}
\right ]
\left [
\begin{array}{c}
\Delta u \\
\Delta p \\
\Delta \phi \\
\Delta \lambda_1 \\
\Delta \lambda_2 \\
\Delta \lambda_3 \\
\end{array}
\right ]
=
\left [
\begin{array}{c}
R_{u}\\
R_{p}\\
R_{\phi}\\
R_{1}\\
R_{2}\\
R_{3}\\
\end{array}
\right ].
\label{eq:orig_mat_struct}
\end{align}
\Marco{where $A$ represents the discretized etc... $G$ represents the discretized gradient operator etc...}\\
\Marco{dimensions matrices?}
\Ido{I don't see how this will clarify things? The matrices are the appropriate jacobians of the formulation. }
In the matrix, presented in block form, the capital letters represent sparse matrices while the small letters denote full vectors.
Here $A$ is the jacobian of the unsteady convection-diffusion part of the momentum equations and $G$ and $G^T$ 
are the discrete gradient and divergence matrices, respectively.
The jacobian $C$ stands for the SUPG formulation of the level-set convection equation, 
while the jacobians $B_1$ and $B_2$ represent the two-way coupling between interface convection and momentum equation. 
The vectors $r_i$, $s_i$ and $t_i$ are the result of the enforcement of the constraints.
The global matrix has a non-symmetric structure due to the absence of the Lagrange multiplier terms in the momentum equations.
Lastly, the right-hand side vector is composed of the residuals of corresponding equations.

Depending on implementation the structure of the matrix can be inconvenient, 
in particular when dealing with a parallel MPI-based solver infrastructure.
Therefore a plain Newton solver is not the solution strategy adopted here, 
instead an alternative  solver strategy analogous to the one discussed in subsection \ref{sec:lm} is used.

\subsubsection{The quasi-Newton solver}
Note that a large portion of the nonlinear character of the problem originates from the additional scalar constraints.
These equations need to be solved to a tight tolerance-level in order to achieve the required conservation behavior.
It is therefore beneficial to decouple the part linked to the constraints from the global problem and solve it to a tight tolerance without incurring high computational costs.
To this purpose we adopt the strategy presented in subsection \ref{sec:lm}.
This results in a convenient matrix structure and allows the  
nonlinear constraints to be solved to a tolerance independent of the tolerance-level of the global problem.
The constraints can be solved with machine precision with minimal overhead.  

Consider the variational formulation (\ref{eq:ce_form}).
Note that each of the Lagrange Multipliers $\lambda_i$ perturbs the level-set solution $\phi$ with a global function, i.e. we write:
\begin {align}
\phi  = \phi_0 + \sum_{i=1}^3 \lambda_i \phi_i,
\end {align}
where $\phi_0$ is the unperturbed solution and $\phi_i~(i=1,2,3)$ are the global perturbations associated with each of the Lagrange multipliers.

By solving the perturbations $\phi_i$ instead of the Lagrange Multipliers $\lambda_i$, the original conservative formulation (\ref{eq:cons_form}) is augmented with:
\\

\textit{Find $ \phi_i \in \WW  $  such that for all $w  \in \WW$,}
\begin{align}\label{eq:pert}
\left (\psi , \frac{\phi_i^{n+1} - \phi_i^{n}}{\Delta t} \right ) 
+ (  \psi ,  \bu^{n+1/2}  \cdot \nabla \phi_i^{n+1/2} )  \nn\\
+\left (  \tau \bu^{n+1/2}  \cdot \nabla \psi ,  \frac{\phi_i^{n+1} - \phi_i^{n} }{\Delta t} + \bu^{n+1/2}  \cdot \nabla \phi_i^{n+1/2} \right ) 
&=  \delta h_i, \quad \quad \text{for} \quad i=1,2,3.
\end{align}
Here $\delta h_i$ are the variations of the constraints with respect to $\phi$ given in (\ref{eq:const_vars 1})-(\ref{eq:const_vars 3}). 
This converts the original matrix structure (\ref{eq:orig_mat_struct}) into to a more standard form:
\begin {align}
\left [
\begin{array}{cccccc}
A      &G &B_1 & 0 &0 &0\\
G^T & 0& 0     & 0 &0 &0\\
B_2&0&C& 0 &0 &0\\
0&0&0& C &0 &0\\
0&0&0& 0 &C &0\\
0&0&0& 0 &0 &C\\
\end{array}
\right ]
\left [
\begin{array}{c}
\Delta u \\
\Delta p \\
\Delta \phi \\
\Delta\phi_1 \\
\Delta\phi_2 \\
\Delta\phi_3 \\
\end{array}
\right ]
=
\left [
\begin{array}{c}
R_{u}\\
R_{p}\\
R_{\phi}\\
R_{\phi_1}\\
R_{\phi_2}\\
R_{\phi_3}\\
\end{array}
\right ],
\label{eq:new_mat_struct}
\end{align}
which only consists of sparse matrix blocks, where $R_{\phi_i}$ are the residuals of (\ref{eq:pert}). 
We solve this system of equations using a standard flexible GMRES with additive Schwartz preconditioning, as provided by Petsc 
\cite{PETSC,petsc-efficient}. The Lagrange Multipliers
 are determined via: 
\begin{align}
h_j \left(\phi_0 + \sum_{i=1,2,3} \lambda_i \phi_i\right) =0,
\label{eq:nl_constraint}
\end{align}
where $h_j$ represents the constraints given in (\ref{eq:disc_mass_constr}), (\ref{eq:disc_kin_eng_constr}) and (\ref{eq:disc_pot_eng_constr}).
This is a nonlinear system of three equations with three unknowns. This is efficiently solved using 
the Newton method at each global iteration. This results in a nested nonlinear iteration loop.
Note that this does not pose any computational problems since the nested problem is very small and thus easily solved.
This subsolver only iterates over the constraints. This results in a tight enforcement of these constraints without the need to
solve the entire coupled problem (\ref{eq:orig_mat_struct}). 

A detailed step-by-step description of the algorithm is given in  \ref{app:algo}.

\section{Numerical results}
\label{sec:results}

In this section we test the performance of proposed energy-corrected formulation (\ref{eq:ce_form}) on a dambreak problem. 
To investigate the importance of the correctness of the kinetic and potential energy evolution, we also carry out simulations without the corresponding constraints. 
In these simulations only the proposed mass correction (\ref{eq:mass_cont}) is active. 
We refer to this conservative method with correct mass behavior as the conservative formulation in the following. 
\added[id=Rev2.1]{To benchmark the numerical results we also employ a formulation in the convective form. 
This is because most two-fluid simulations based on the level-set methodology  
are performed with this form, see for instance \cite{Akin05a,LoYaOn06,ElCo07,LERC10,CrCeTe07,CrCeTe07b}}. 
We employ a convective formulation with correct mass behavior.
This results in a method that is very similar to  previously published \cite{KAFB10,ABKF11,AkBaBeFaKe12,ADKSB12}.
For the precise convective formulation consult \ref{app:conv_form}.

\subsection{Dambreak problem}

A well-known dambreak problem serves as test case for the verification of the energy evolution of the presented methods.
The employed setup closely resembles the one of Martin and Moyce \cite{MaMo52,KoOk96}.
In this problem a column of water of size $14.6~\text{cm}$ by $29.2~\text{cm}$ collapses due to the gravitational force.
The computational domain is of size $58.4~\text{cm}$ by $35.04~\text{cm}$ and filled for the remainder with air.
For the densities of the fluids we use  $\rho_0  = 1.00~\text{kg/m$^3$}$  and $\rho_1=1000~\text{kg/m$^3$}$, which is similar to the physical values for air and water.
The viscosity is set to $\mu_0=\mu_1=2.0~\text{kg/(m.s)}$ for both fluids.
This is significantly higher than the physical value for the dynamic viscosity.
These higher values are chosen in order to avoid instabilities that would otherwise occur in  the Galerkin discretization. 
Stabilized formulations are avoided on purpose, as correct energy evolution of stabilized formulations is a problem by itself.
Consult \cite{EiAk17i,EiAk17ii} for further elaboration and solution strategies.

All computations are done on uniform Cartesian meshes with mostly linear shape functions. 
Each velocity component is discretized with a $C_1$-quadratic shape functions in the appropriate direction.
This deviation is essential to arrive at a stable velocity-pressure pair that results in solenoidal solutions.
Furthermore, the computations are performed with no-penetration boundary conditions on all the surfaces.
 
All computations are done with the Crank-Nicolson time-integration. The time-step is adjusted with a simple 
proportional controller \cite{VaCaCo05},
\begin{align} 
\Delta t^{n+1} = \left (\frac{{\rm CFL_{target}}}{{\rm CFL}^n}\right)^{K_p} \Delta t^n
\end{align}  
where the {\rm CFL}-number is defined as,
\begin{align} 
{\rm CFL}^n = \Delta t^n \max_{\bx \in \Omega}  \sqrt{\bu^n\cdot \B{G}\bu^n}
\end{align}  
and ${\rm CFL_{target}}$ is the target value and $K_p$ is the proportional gain. Here, both values are set to $0.75$.
\begin{figure}[!ht]
    \begin{center}
    \begin{subfigure}[b]{0.325\textwidth}
        \center
        \includegraphics[width=0.9\textwidth]{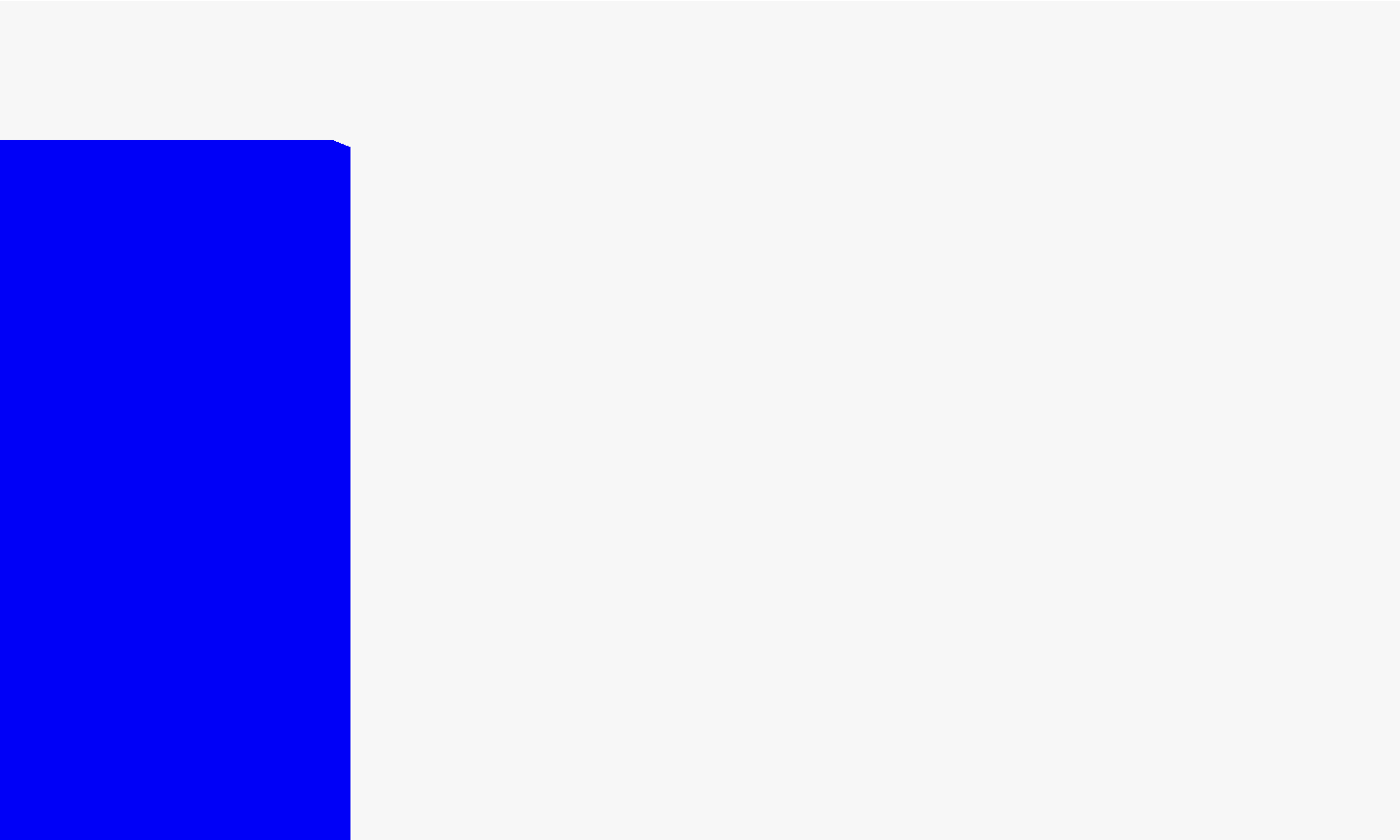}
        \caption{$t=0.0\,s$.}
    \end{subfigure}
    \begin{subfigure}[b]{0.325\textwidth}
        \center
        \includegraphics[width=0.9\textwidth]{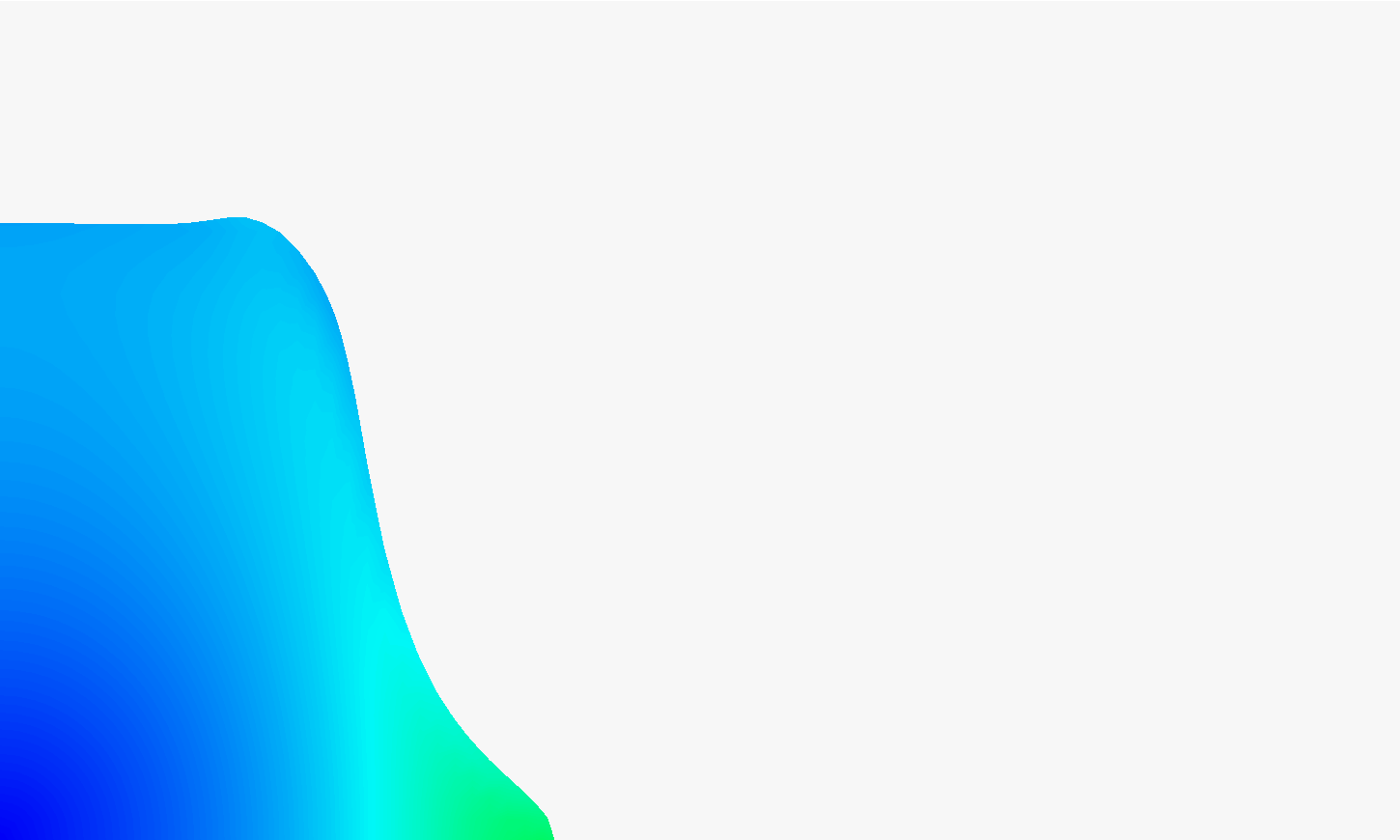}
        \caption{$t=0.1\,s$.}
    \end{subfigure}
     \begin{subfigure}[b]{0.325\textwidth}
       \center
       \includegraphics[width=0.9\textwidth]{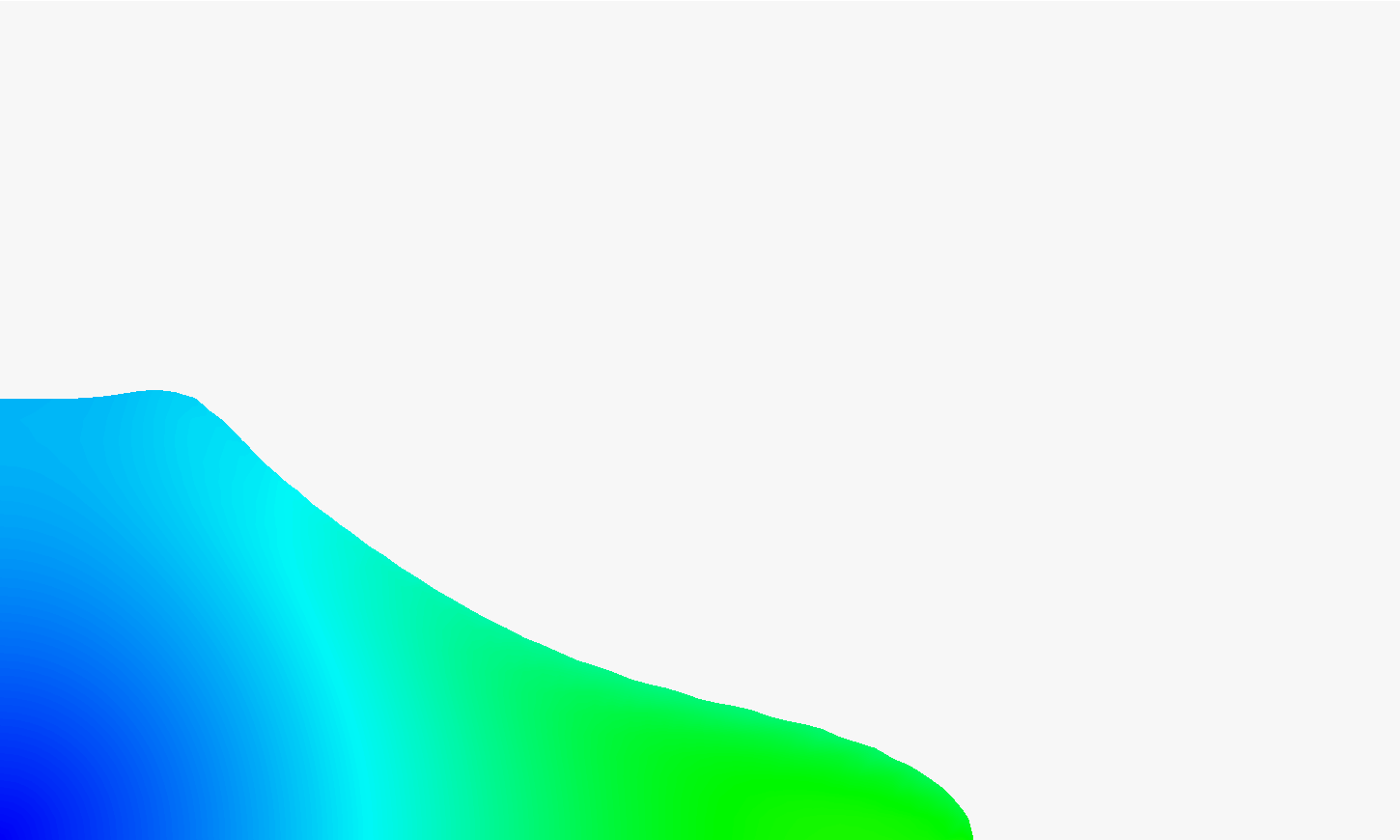}
        \caption{$t=0.2\,s$.}
    \end{subfigure} 
    \\~\\
    \begin{subfigure}[b]{0.325\textwidth}
        \center
        \includegraphics[width=0.9\textwidth]{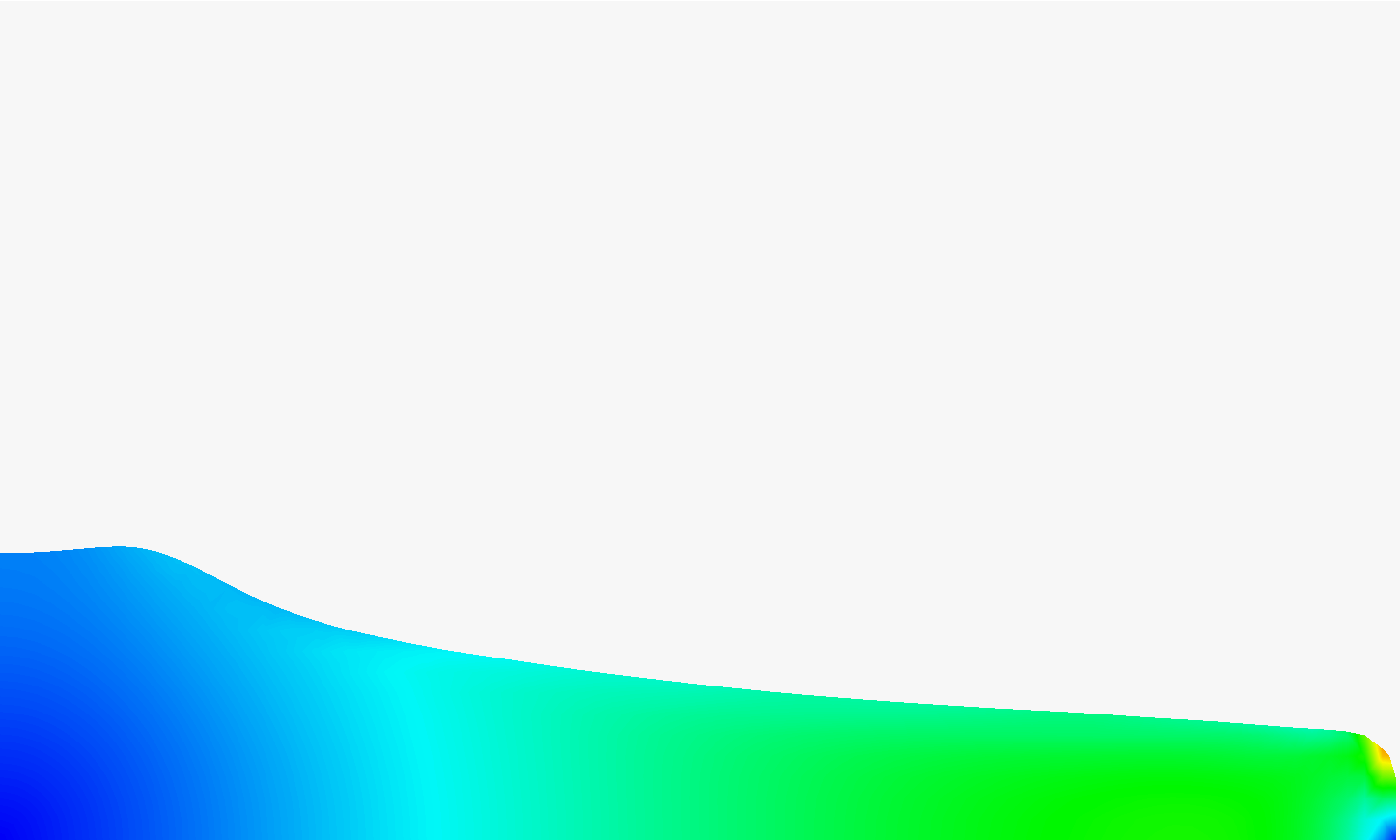}
        \caption{$t=0.3\,s$.}
    \end{subfigure}
    \begin{subfigure}[b]{0.325\textwidth}
        \center
        \includegraphics[width=0.9\textwidth]{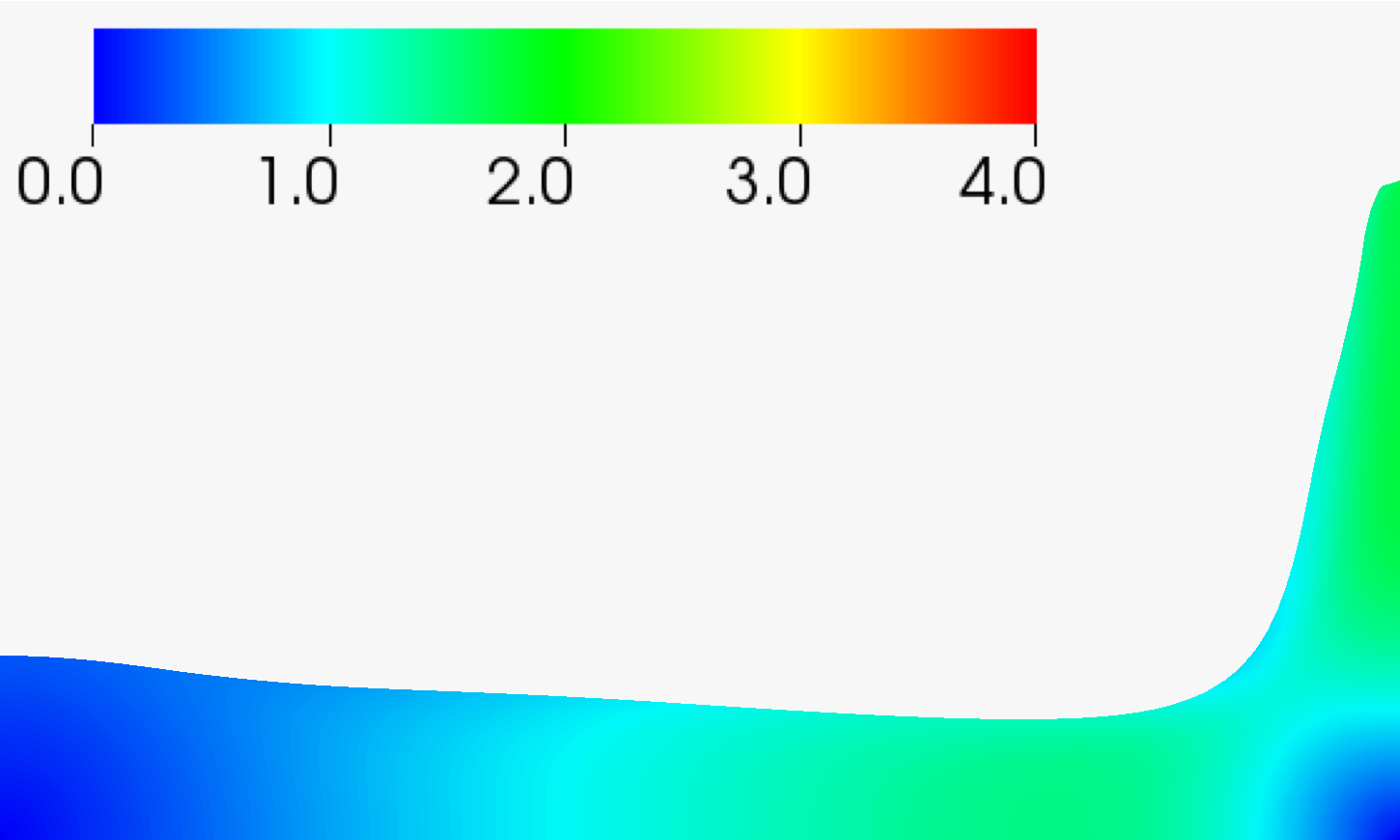}
        \caption{$t=0.4\,s$.}
    \end{subfigure}
    \begin{subfigure}[b]{0.325\textwidth}
        \center
       \includegraphics[width=0.9\textwidth]{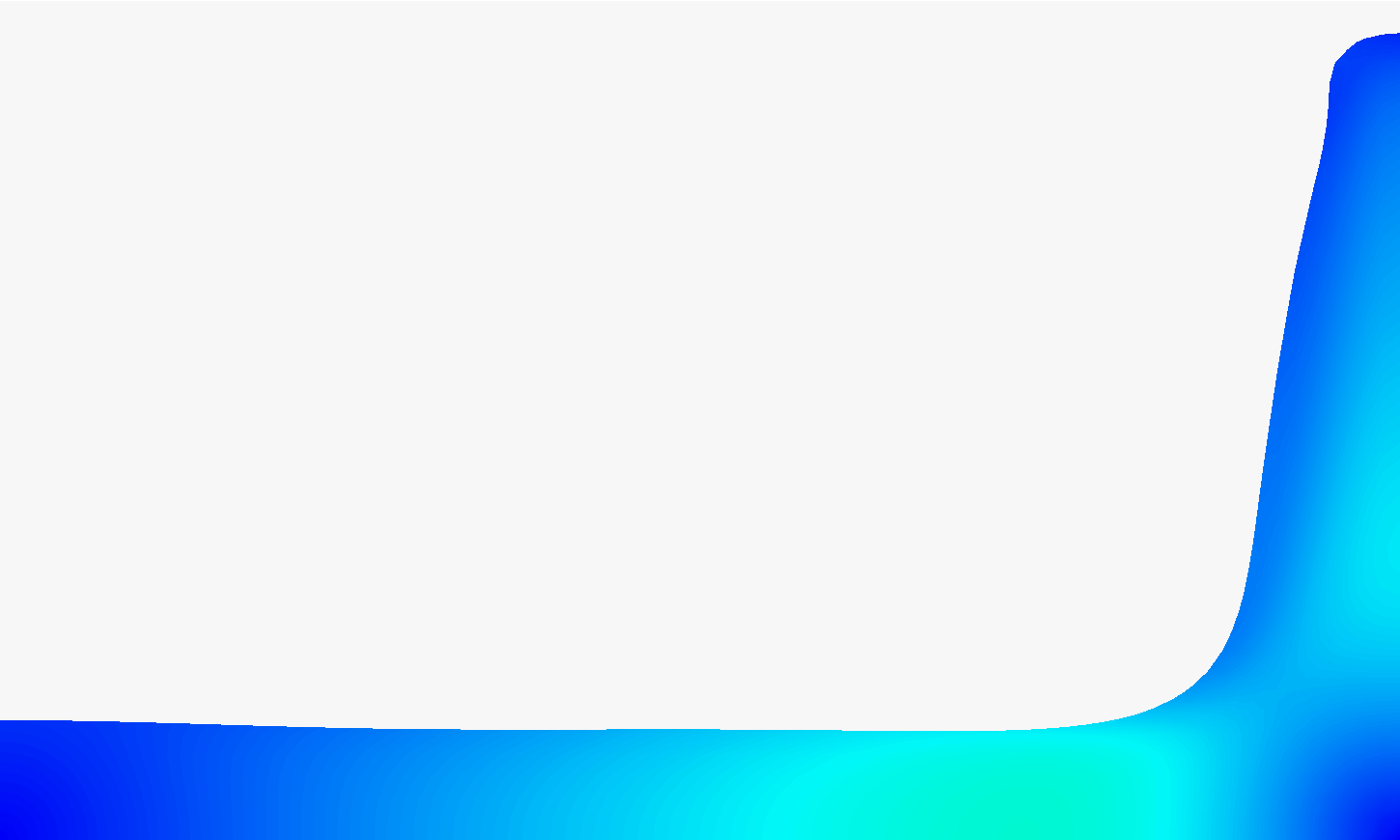}
        \caption{$t=0.5\,s$.}
    \end{subfigure}
    \\~\\
    \begin{subfigure}[b]{0.325\textwidth}
        \center
        \includegraphics[width=0.9\textwidth]{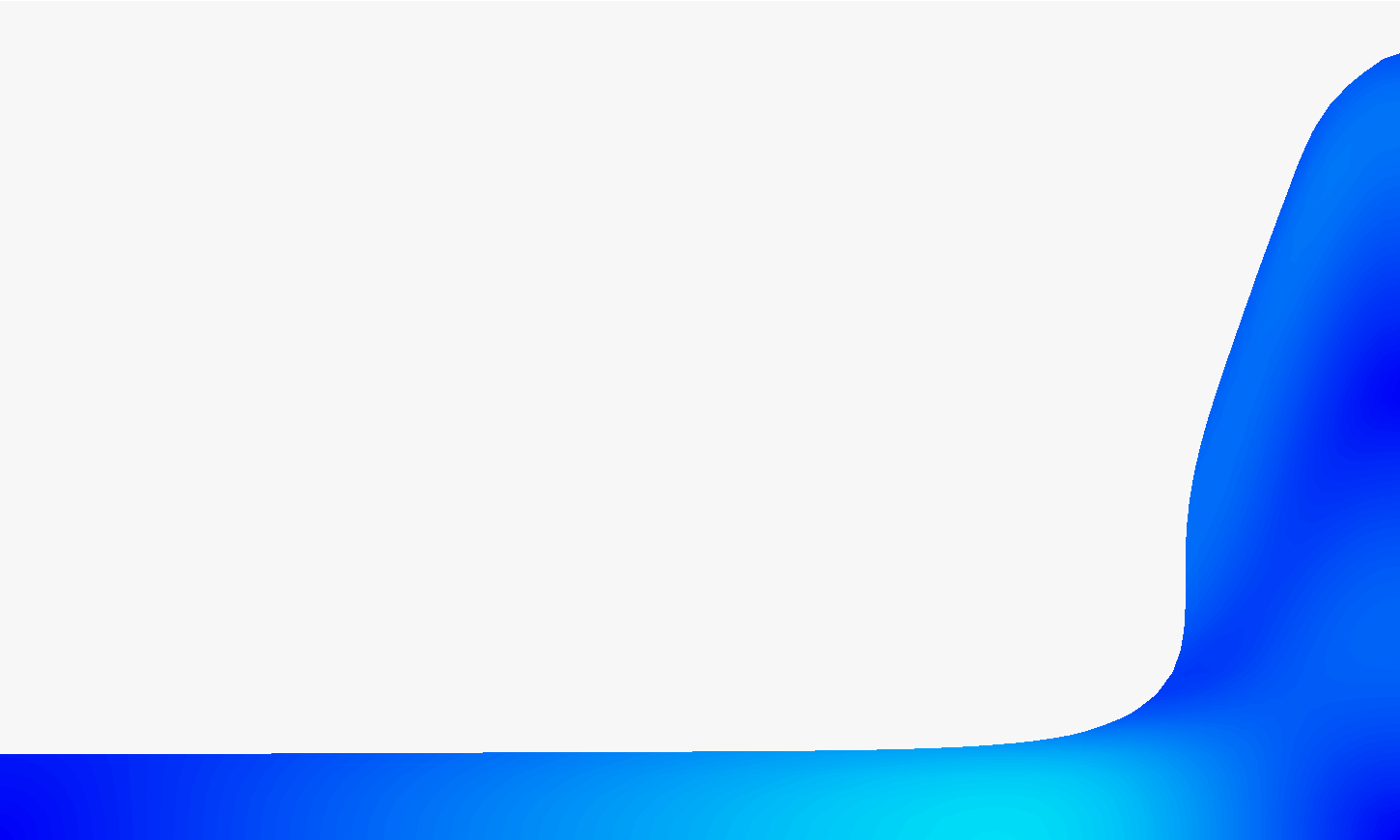}
        \caption{$t=0.6\,s$.}
    \end{subfigure}
    \begin{subfigure}[b]{0.325\textwidth}
        \center
        \includegraphics[width=0.9\textwidth]{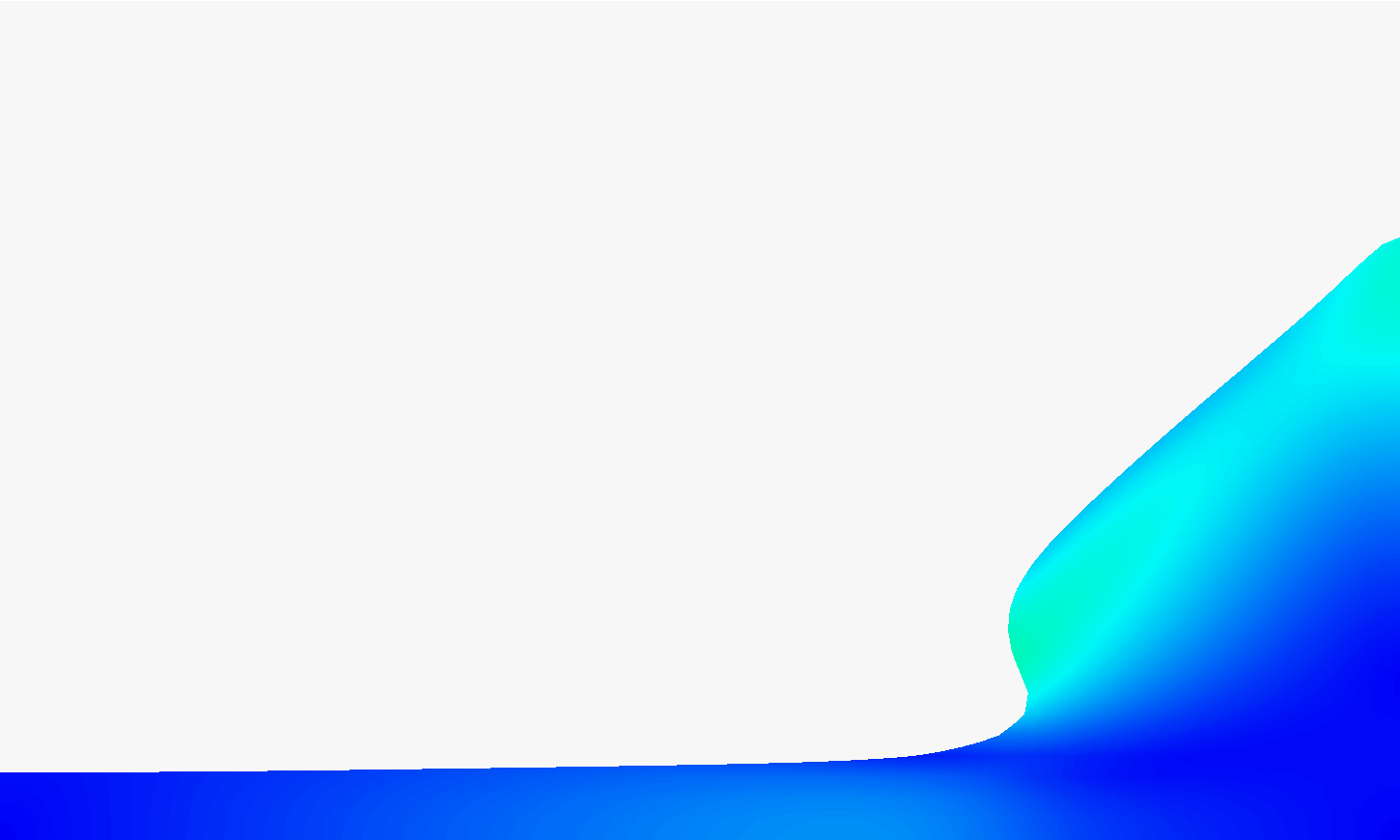}
        \caption{$t=0.7\,s$.}
    \end{subfigure}
    \begin{subfigure}[b]{0.325\textwidth}
        \center
        \includegraphics[width=0.9\textwidth]{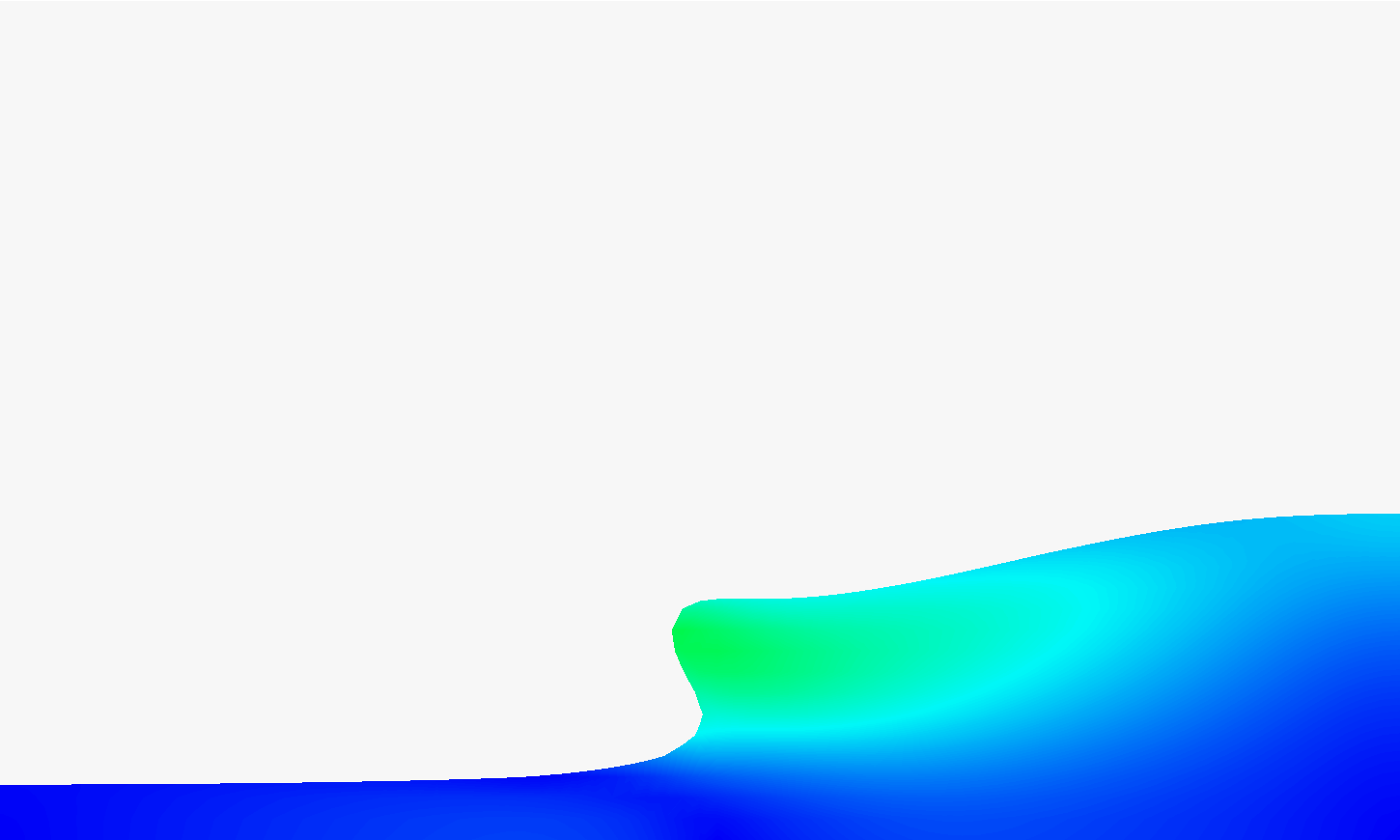}
        \caption{$t=0.8\,s$.}
    \end{subfigure}
    \caption{Snapshots of the time evolution of the water column in the dambreak problem. 
    The colors represent the velocity magnitude given in m/s.}\label{fig:snap}
    \end{center}
\end{figure}

Figure \ref{fig:snap} shows 9 snapshots of the solution profiles.
The water column collapses, subsequently moves rightwards, runs up the right wall and moves back leftwards.
All the snapshots display smooth solutions.
This is linked to (i) the relatively high viscosity employed, (ii) the smooth NURBS basis functions and (iii) the novel level-set formulation.

\begin{figure}[ht!]
    \begin{center}
    \includegraphics[width=0.48\textwidth]{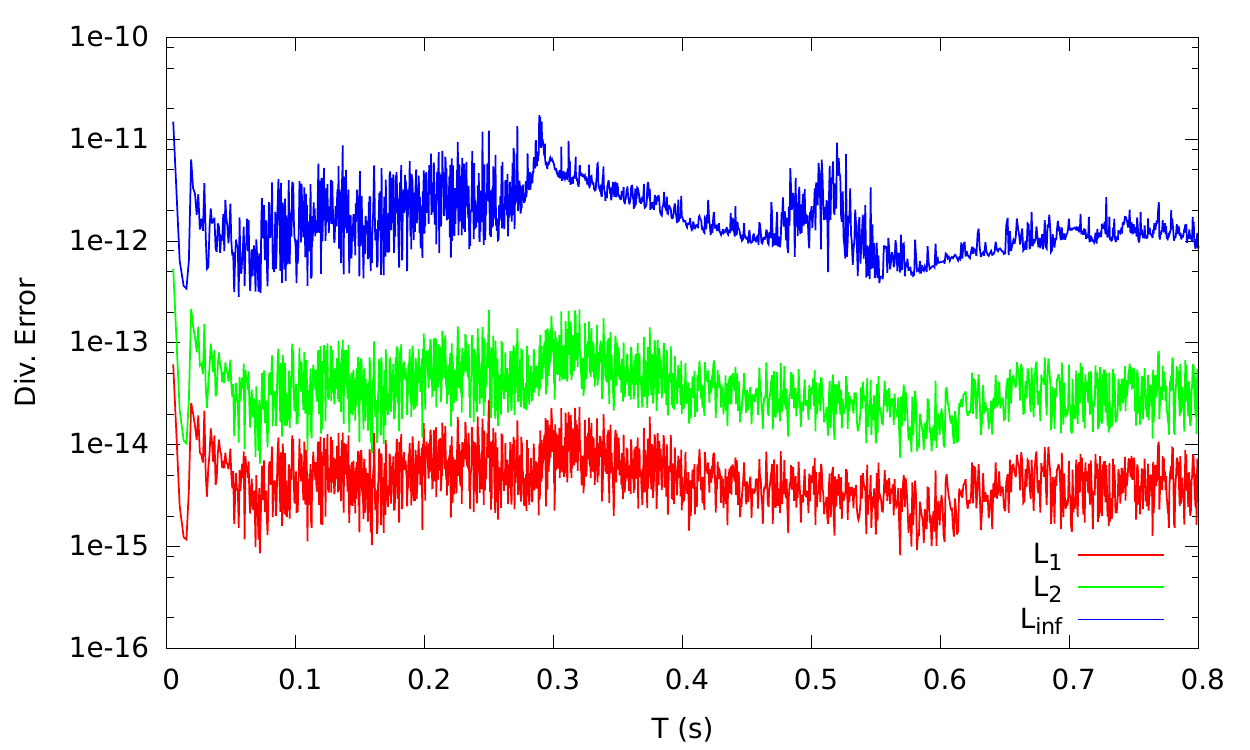}
    \caption{Time trace of the divergence error on a $160$x$80$ mesh.}\label{fig:div}
    \end{center}
\end{figure}
Figure \ref{fig:div} shows a typical the time trace of the divergence error.
This  plot confirms the divergence free character of the solution as discussed in section \ref{sec:div_free}.
The $L_1(\Omega)$-, $L_2(\Omega)$- and $L_{\rm inf}(\Omega)$- errors are all below $10^{-10}$, 
which effectively means that these are almost zero \added[id=Rev2.1]{with respect}  to machine precision. 
Note that the $L_{\rm inf}$-error is determined \added[id=Rev2.1]{by sampling the error at the integration points.
Given that the divergence is in the pressure space, this implies it must be almost zero up to machine precision everywhere.
}

\subsection{Examination of the energy behavior}
\begin{figure}[h!]
    \begin{center}
    \includegraphics[width=0.48\textwidth]{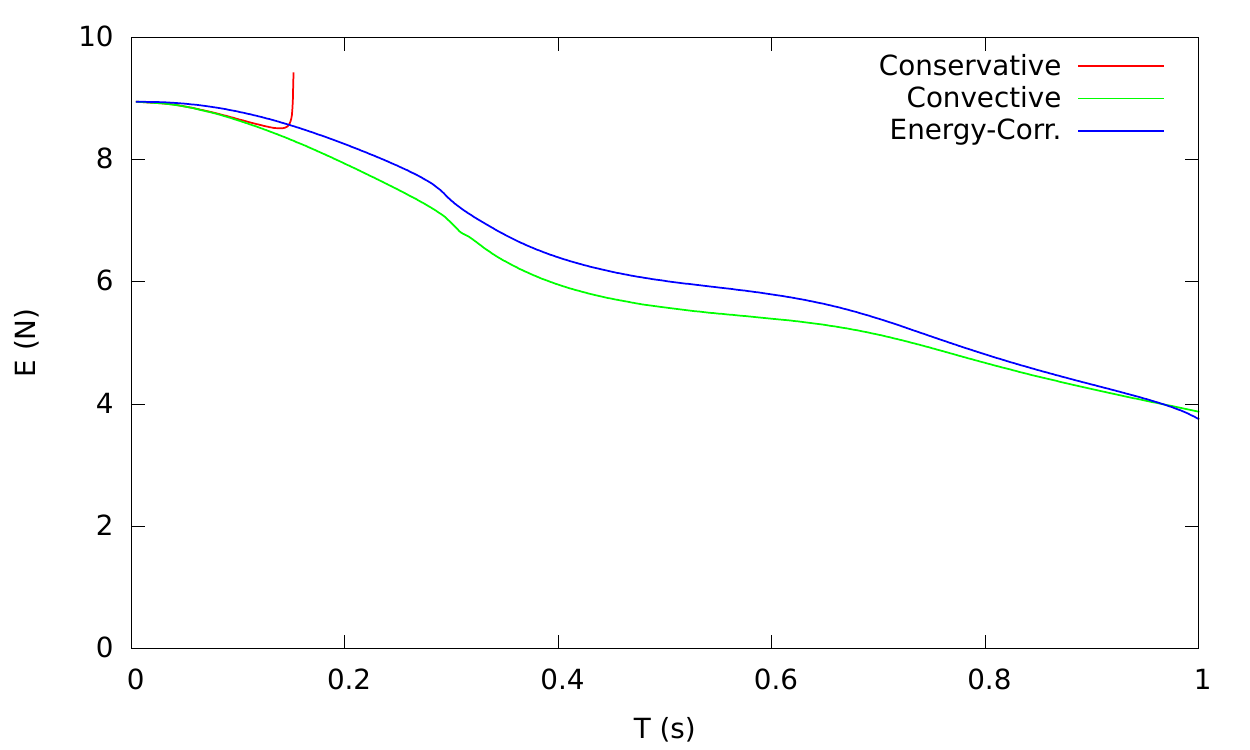}
    \caption{Total energy decay  for the different methods on a $40$x$20$ mesh.}\label{fig:m40}
    \end{center}
\end{figure}

Here we test the energy behavior of the various formulations. 
Figure \ref{fig:m40} visualizes the energy evolution of the different formulations.
During the simulation no energy is supplied to the system, i.e. the total energy needs to decay monotonically.
This is not the case for the conservative formulation which shows an exponential growth of energy.
The result is an unstable and diverging computation.
We do not further consider this approach.
In contrast, the monolithic convective and monolithic energy-corrected 
formulations do result in stable simulations with monotonic energy decay.
The energy decay of the convective formulation is somewhat misleading as it is partly the consequence of the high viscosity in the presented test case. 
Note that the conservative and energy-corrected formulation show a small difference in energy behavior.
We now look into this gap.

\begin{figure}[!ht]
    \begin{center}
    \begin{subfigure}[b]{0.48\textwidth}
    \includegraphics[width=\textwidth]{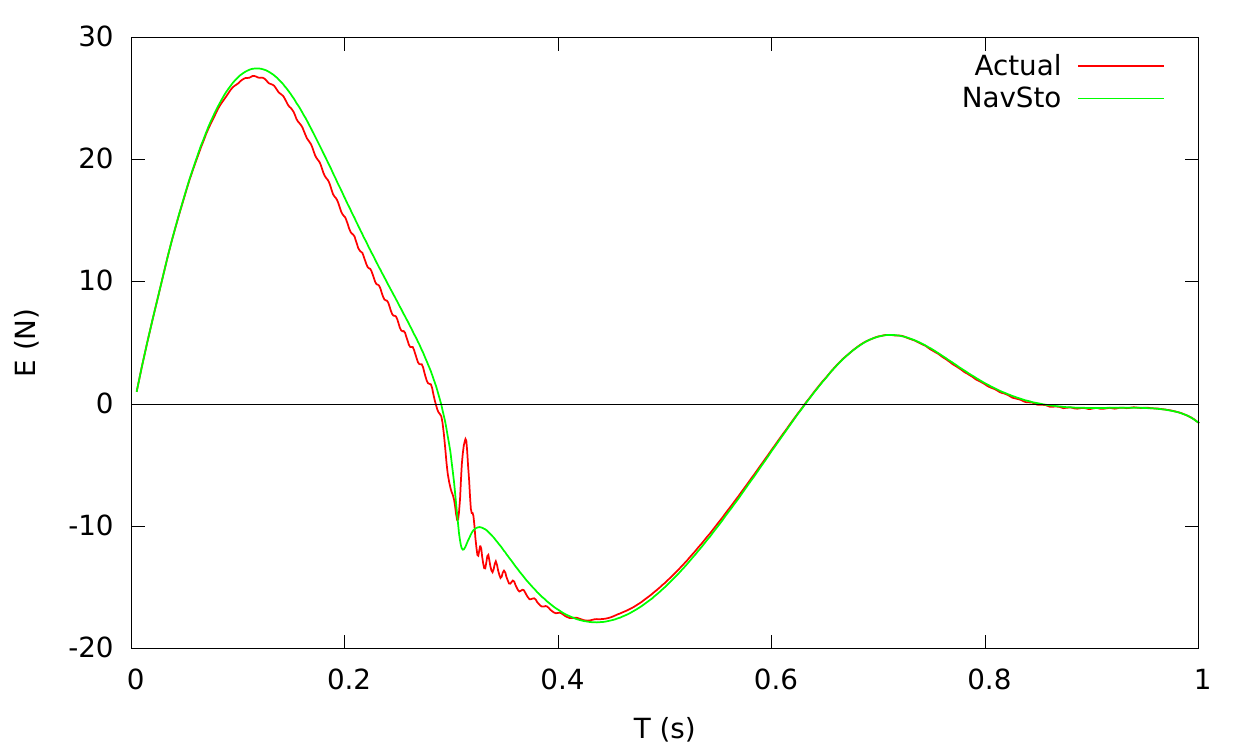}
        \caption{Rate of change of kinetic energy.}
    \end{subfigure}
    \begin{subfigure}[b]{0.48\textwidth}
    \includegraphics[width=\textwidth]{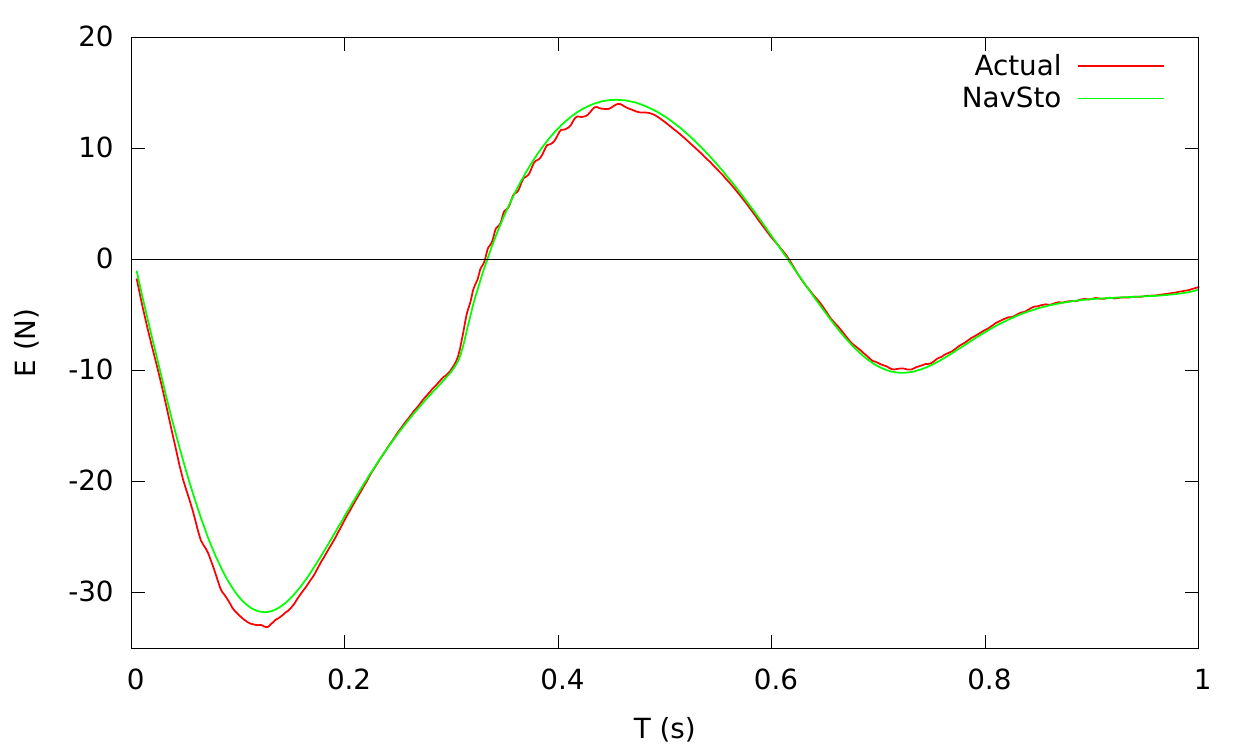}
        \caption{Rate of change of potential energy.}
    \end{subfigure}
    \caption{Rate of change of energies for the convective formulation.}\label{fig:eng_conv}
    \end{center}
\end{figure}

Figure \ref{fig:eng_conv} shows the rate of change in kinetic and potential energy for the monolithic convective formulation.
This rate of change is computed in two different ways.
The first one is the actual difference between the energies in two consecutive time steps:
\begin {align}
\dfrac{{\rm d}}{{\rm d}t}E_{pot} =& \Delta t ^{-1}\left (\int_\Omega \rho^{n+1} \B{g}\cdot \bx {\rm d}\Omega -\int_\Omega \rho^n \B{g}\cdot \bx {\rm d}\Omega \right), \\
\dfrac{{\rm d}}{{\rm d}t}E_{kin} =& \Delta t ^{-1}\left (\int_\Omega \onehalf \rho^{n+1} \bu^{n+1}\cdot \bu^{n+1} {\rm d}\Omega -\int_\Omega \onehalf \rho^n \bu^n\cdot \bu^n {\rm d}\Omega \right), 
\end{align}
which is referred to as actual. The second rate of change is a direct result of the weak formulation of the Navier-Stokes equations (\ref{eq: convective form}).
In a computation with correct energy behavior these approaches should provide the same results.
\\
Figure \ref{fig:eng_conv} displays a mismatch of the two rates of energy.
Globally they follow the same trend but there are some clear deviations, particularly around $t=0.3~\text{s}$ 
when the water hits the right wall. 
Furthermore, the actual rate of energy change shows some wiggles.

\begin{figure}[!ht]
    \begin{center}
    \begin{subfigure}[b]{0.48\textwidth}
    \includegraphics[width=\textwidth]{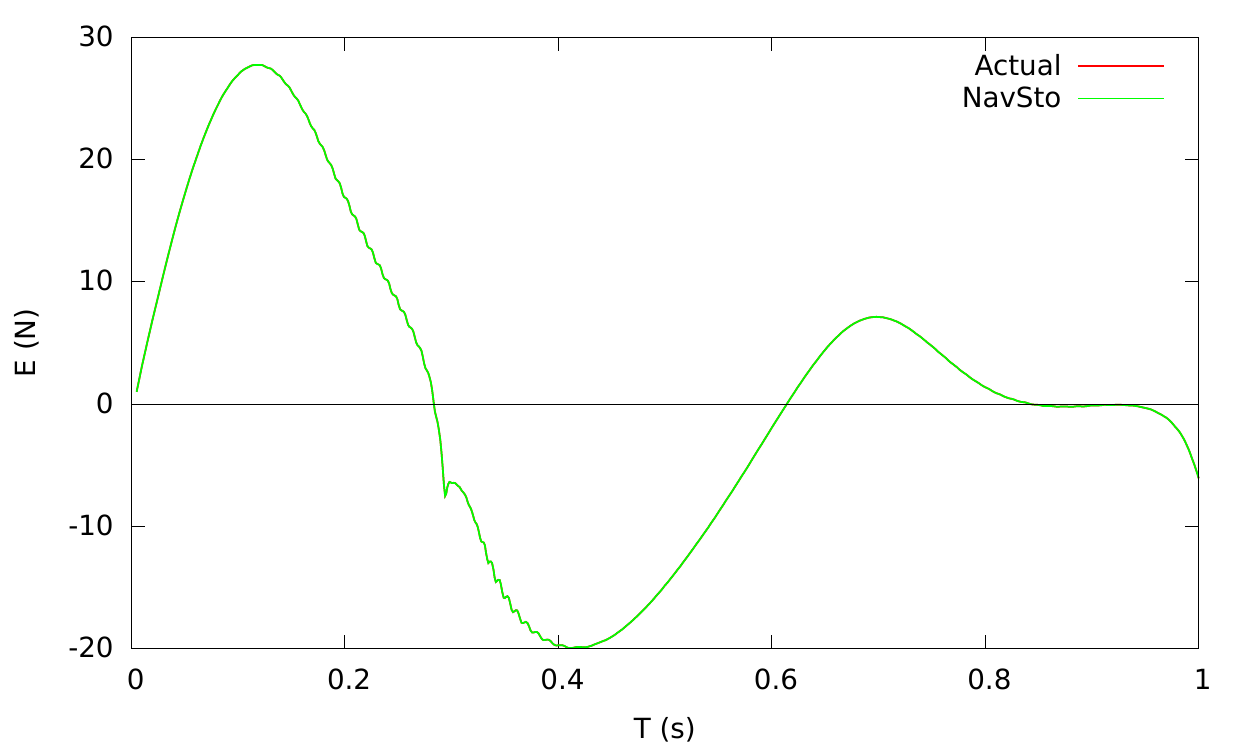}
        \caption{Rate of change of kinetic energy.}
    \end{subfigure}
    \begin{subfigure}[b]{0.48\textwidth}
    \includegraphics[width=\textwidth]{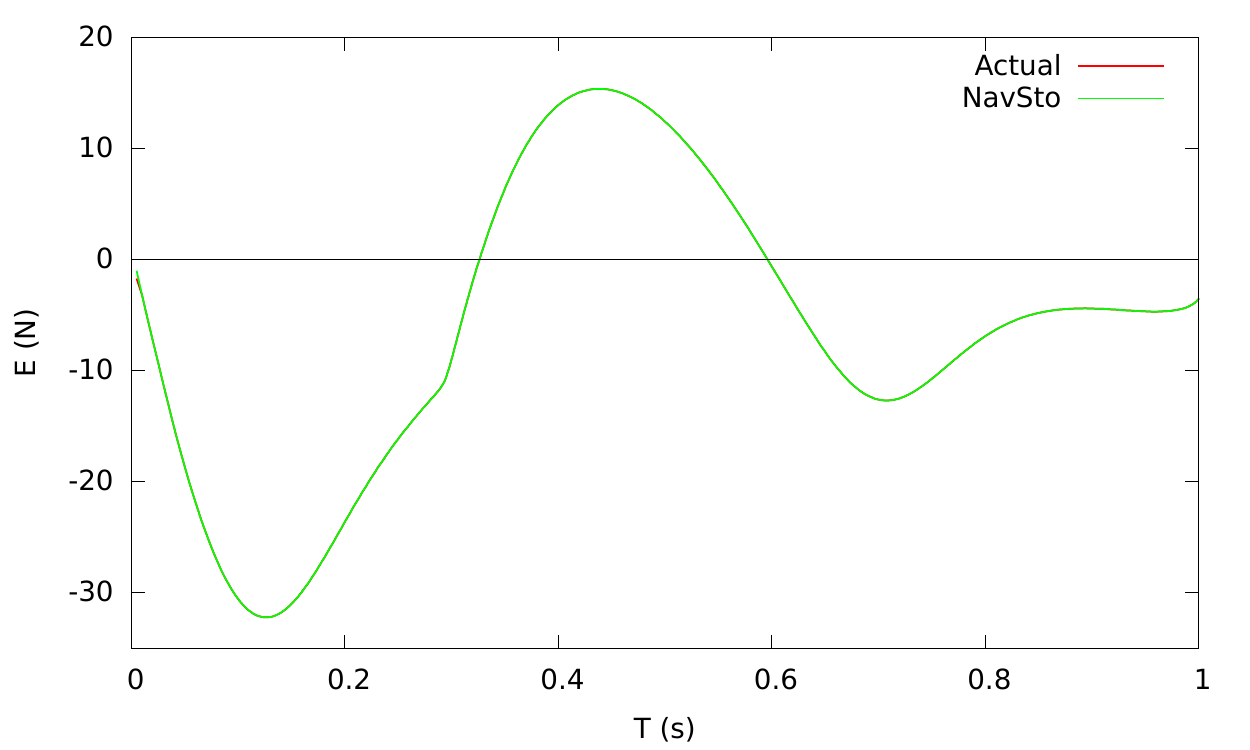}
        \caption{Rate of change of potential energy.}
    \end{subfigure}
    \caption{ Rate of change of energies for the energy--corrected formulation.}\label{fig:eng_lm}
    \end{center}
\end{figure}

Next, in Figure \ref{fig:eng_lm} we depict the same rates for the energy-corrected formulation.
Here the two rates of change are exactly the same. 
The kinetic and potential energy constraints force the interface to evolve such that
the correct global energy behavior is obtained. 
Note that the wiggles in the rate of change of the potential 
energy have vanished, whereas the wiggles in kinetic energy evolution rate are significantly attenuated.
This clearly shows the importance of correct energy behavior for the overall quality of the obtained solutions.

\begin{figure}[!ht]
    \begin{center}
    \begin{subfigure}[b]{0.48\textwidth}
    \includegraphics[width=\textwidth]{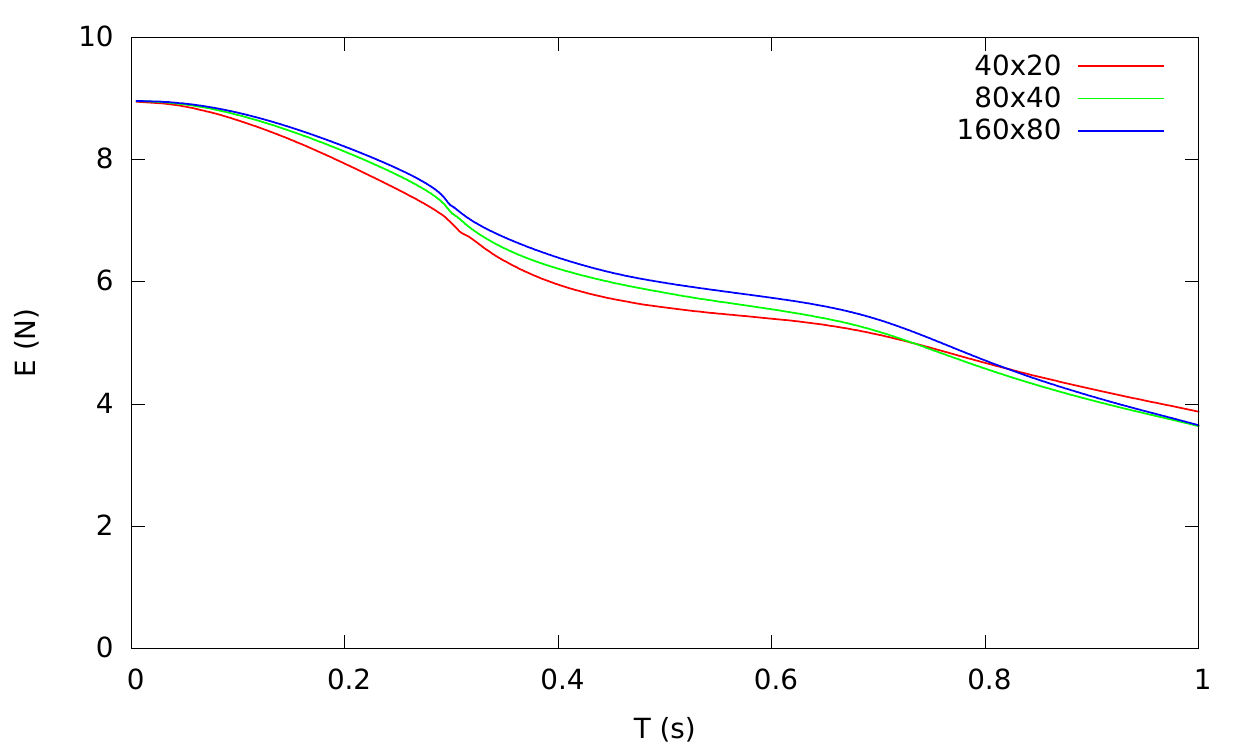}
        \caption{Convective formulation.}
    \end{subfigure}
    \begin{subfigure}[b]{0.48\textwidth}
    \includegraphics[width=\textwidth]{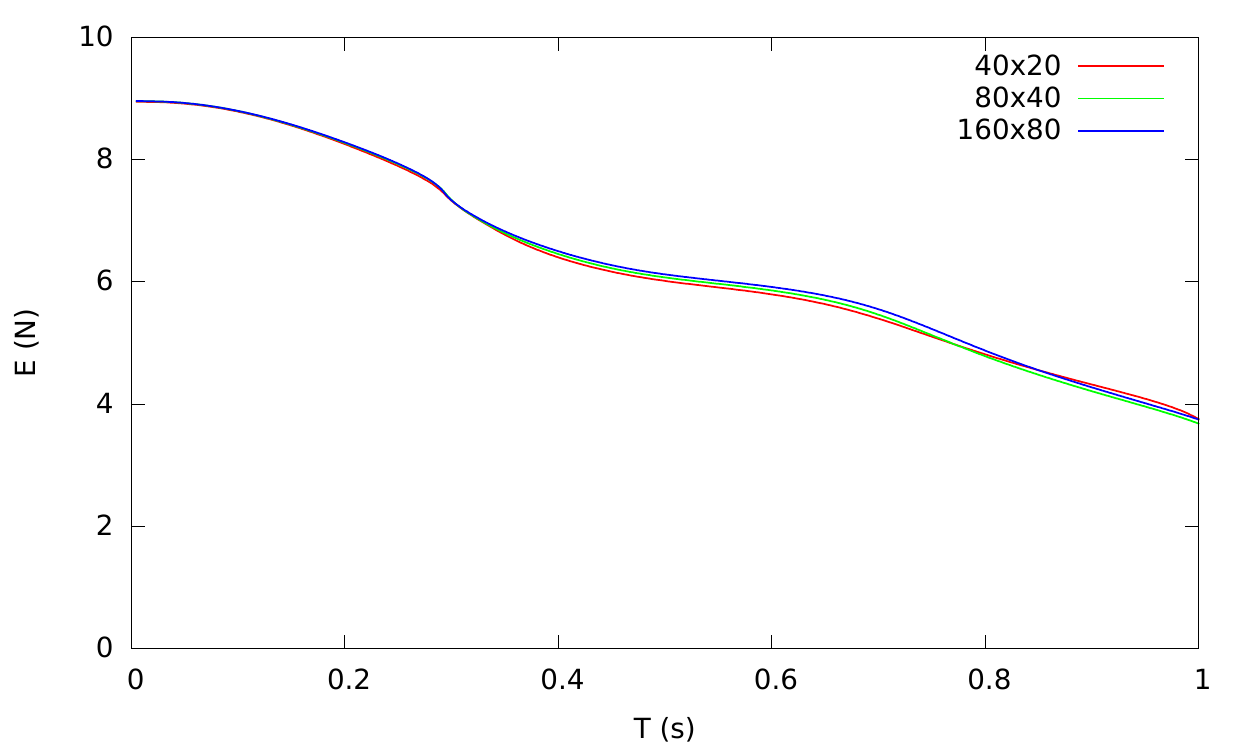}
        \caption{Energy--corrected formulation.}
    \end{subfigure}
    \caption{Convergence of the energy decay.}\label{fig:conv_lm}
    \end{center}
\end{figure}

\begin{figure}[!ht]
    \begin{center}
    \begin{subfigure}[b]{0.48\textwidth}
    \includegraphics[width=\textwidth]{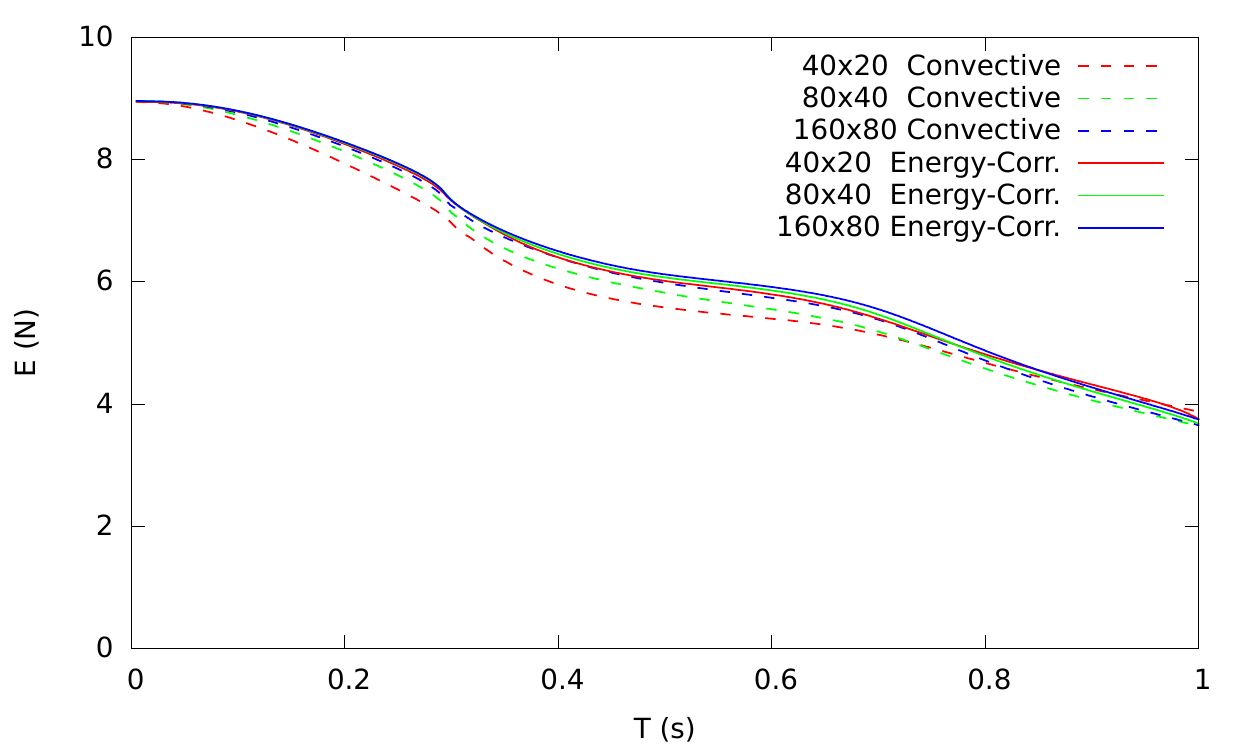}
    \caption{Convective and energy--corrected formulation.}\label{fig:conv_lm2}
    \end{subfigure}
    \begin{subfigure}[b]{0.48\textwidth}
    \includegraphics[width=\textwidth]{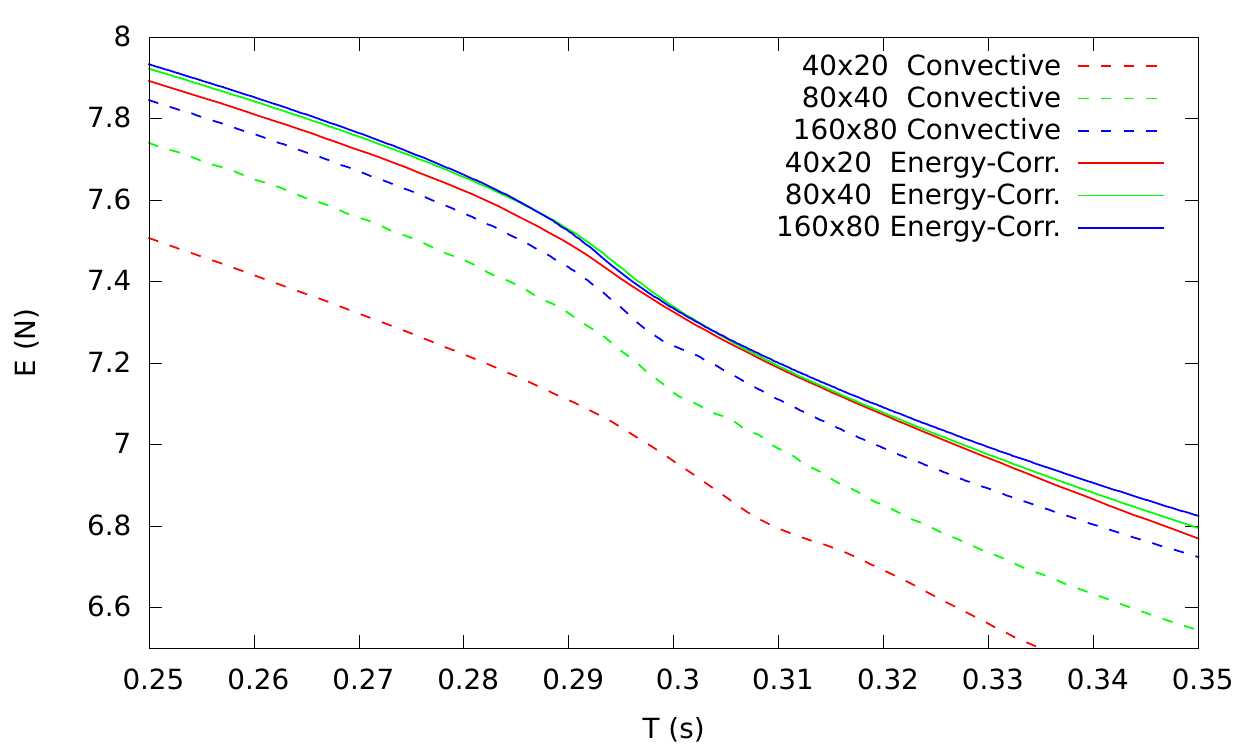}
    \caption{Zoom of same plot. }\label{fig:conv_lm2_zoom}
    \end{subfigure}
    \caption{Convergence of the energy decay.}\label{fig:conv_lm_o}
    \end{center}    
\end{figure}

Lastly, we focus on the convergence characteristics of the methods. Figure \ref{fig:conv_lm} shows the convergence plots of the global energy evolution. 
Both methods clearly converge. 
However, the energy-corrected method on the coarser meshes already exhibits an excellent agreement with results on the finer meshes. 
This is not the case for the convective method.  
This is also clear from Figure \ref{fig:conv_lm_o} which shows the convergence plots of the two methods in one plot.
The energy evolution of the energy-corrected method on the coarsest mesh ($40$x$20$) is very similar to that of the convective method on the finest mesh ($160$x$80$). 
Hence, the novel approach improves significantly the accuracy.

\section{Conclusion}
\label{sec:conclusion}
In this work we have presented a method with correct energy behavior for the computation of two-fluid flow.
The formulation is of conservative-type and uses the level-set method to describe the interface.

The analysis of the conservation properties (mass, momentum and energy) on a continuous level reveals that the correct evolution of interface is of critical importance.
In contrast to the continuous form, the standard discrete formulation does not guarantee these conservation properties.
This is linked to both the spatial and temporal integration of the interface evolution.

The proposed methodology rectifies these discrepancies by explicitly enforcing correct mass, kinetic energy and potential energy behavior in the formulation.
These constraints are enforced via a Lagrange multiplier construction in the interface evolution equation.
The level-set convection equation itself is stabilized with a standard SUPG approach.
Furthermore, the approach is presented in the isogeometric analysis framework to ensure exact incompressiblity of the velocity fields.
This is a natural feature of the presented method since it reduces the approach to a valid method in the single-fluid case.

The implementation employs a quasi-newton approach to solve the nonlinear system.
This approach partially disconnects the constraints from the rest of the global problem.
It leads to a favorable matrix structure, isolates the `most' nonlinear part of the formulation, and 
allows strict enforcement of the constraints without incurring too much computational overhead.

We have tested the presented methods on a prototype dambreak problem.
The numerical results show that the standard conservative method breaks down whereas the novel methodology shows excellent performance.
A standard convective formulation, serving as reference result, provides reasonable results.
However when looking at the kinetic and potential energy evolution there is mismatched between the actual change in energies 
and those experienced by the discretized Navier-Stokes equations. 
The proposed energy--corrected formulation does not have this mismatch and as such has a guaranteed decay of energy.
Furthermore, the novel method requires a significantly smaller amount of grid points compared to the convective formulation.
These observations indicate the importance of correct energy behavior in two-fluid flow simulations.
We believe that the large accuracy gain of the new method outweighs its additional implementation effort.

The current  formulation is based on a Galerkin formulation and is therefore only suitable for 
low Reynolds number flows. 
Current work concerns the development of a two-fluid stabilized formulation suitable for the computation of high-Reynolds-number flow problems.

\section*{Acknowledgements}
\label{sec:ack}

The authors are grateful to Delft University of Technology
for its support. Visualizations are done using VisIt \cite{VisIt}, provided by Lawrence Livermore National Laboratory.


\appendix
\section{Quasi-Newton Algorithm}
\label{app:algo}
For clarity a step-by-step description of
 the routine sketched in the previous section  is given here.
The routine describes how the solution at a new time step,  $\bu^{n+1}, p^{n+1}, \phi^{n+1}$ is obtained from the solution at the
current time step,  $\bu^n, p^n, \phi^n$. \\

\noindent The algorithm reads:
 
\begin{enumerate}
\item Start: $\bu^n, p^n, \phi^n$.
\item Initialize the solution:
\begin{align}
\bu^{n+1} =& \bu^n, \nn \\
p^{n+1}    =& p^n, \nn \\
\phi_0^{n+1} =& \phi^n.
\end{align}
\item Initialize the perturbations and the Lagrange multipliers:
\begin{align}
\phi^{n+1}_i =0, \quad i=1,2,3\\
\lambda_i =0, \quad i=1,2,3.
\end{align}

\item \label{it:start} Assemble the right-hand side of (\ref{eq:new_mat_struct}) by evaluating the residuals given by (\ref{eq:cons_form}) and (\ref{eq:pert}).
\item Compute the global norm of the residuals:
$\|R\|^2 =  \|R_{u}\|^2  + \|R_{p}\|^2 + \|R_{\phi}\|^2 + \|R_{\phi_1}\|^2 + \|R_{\phi_2}\|^2 + \|R_{\phi_3}\|^2$.
\item Check convergence: if $\|R\| < \epsilon_1 \|R\|_{\text{ref}}$ then go to step \ref{it:finish}.
\item Assemble the matrix given in (\ref{eq:new_mat_struct}) by evaluating the Jacobians of (\ref{eq:cons_form}) and (\ref{eq:pert}).
\item Solve the linear problem given in (\ref{eq:new_mat_struct}).
\item Update the solution:
\begin{align}
\bu^{n+1} =& \bu^{n+1}  + \Delta \bu,     \\    
p^{n+1}    =& p^{n+1}  + \Delta p,         \\   
\phi_0^{n+1} =& \phi_0^{n+1} + \Delta \phi_0, \\
\phi_i^{n+1} =& \phi_i^{n+1} + \Delta \phi_i, \quad i=1,2,3.
\end{align}

\item Solve the nonlinear system (\ref{eq:nl_constraint}):
\begin{enumerate}
\item Assemble the right-hand side by evaluating (\ref{eq:nl_constraint}).
\item Compute norm and check convergence:
\vspace{-2.5mm}
\begin{align}
 \sum_{i=1,2,3} h_i^2 < \epsilon_2^2
\end{align}
\item Assemble the matrix by evaluating the Jacobian of (\ref{eq:nl_constraint}).
\item Solve the linear problem using a direct solver.
\item Update the solution:
\vspace{-2.5mm}
\begin{align}
\lambda_i = \lambda_i + \Delta \lambda_i, \quad i=1,2,3.
\end{align}
\end{enumerate}

\item Go to step \ref{it:start}.

\item \label{it:finish} Update the level-set variables:
\begin{align}
\phi^{n+1} =& \phi_0^{n+1} + \sum_{i=1,2,3} \lambda_i \phi^{n+1}_i,
\\
\phi_i^{n+1}=& 0, \quad i=1,2,3.
\end{align}

\item Finish: $\bu^{n+1}, p^{n+1}, \phi^{n+1}$.
\end{enumerate}

The residual norm of the first iteration is used as the reference residual norm $\|R\|_{\text{ref}}$.
The convergence tolerances are typically taken as $\epsilon_1=10^{-3}$ and $\epsilon_2=10^{-12}$. 
Note that $\epsilon_1$ is a relative tolerance and that $\epsilon_2$ is an absolute tolerance.

\section{Standard convective discretization}
\label{app:conv_form}

Here we present the standard convective discrete formulation which serves as a benchmark method. 
The convective form in strong form follows when applying the incompressibility constraint in the momentum equation:
\begin{subequations}
\begin{align} 
 \partial_t (\rho \bu) +  \rho \bu \cdot \nabla \bu
 + \nabla  p - \nabla \cdot 2\mu \nabla^s \bu &=  \rho \B{g}, \\
    \nabla \cdot \bu &= 0, \\
 \partial_t \rho + \bu \cdot \nabla \rho &= 0.
\end{align}
\end{subequations}
Using this strong form the standard discrete formulation in convective form with mass conservation reads:\\


 \textit{Find $\mathbf{u}^{n+1} \in \mathcal{U} ,~p^{n+1} \in \mathcal{P} ,~\phi^{n+1} \in H^1(\Omega)$ 
and $\lambda_1 \in   \mathbb{R}$ such that for all}
\textit{$~\mathbf{w} \in \mathcal{U},~q \in \mathcal{P} ,~\psi \in H^1(\Omega)$,}

\begin{subequations}\label{eq: convective form}
\begin{align}
\left (\bw,  \rho^{n+1/2} \frac{\bu^{n+1} -\bu^n}{\Delta t}\right )
 +( \bw, \rho^{n+1/2} \bu^{n+1/2}\nabla \bu^{n+1/2})  
 - (\nabla \cdot \bw, p^{n+1}) &\nn\\ 
 + (\nabla \bw, 2\mu \nabla^s \bu^{n+1/2}) &= (\bw, \rho^{n+1/2} \B{g}), \\
 ( q, \nabla \cdot \bu^{n+1/2}) &= 0, \\
\left (\psi, \frac{\phi^{n+1} - \phi^{n} }{\Delta t}+   \bu^{n+1/2}  \cdot \nabla \phi^{n+1/2} \right)
 \qquad\qquad\qquad\qquad\qquad\qquad\qquad \nn\\
+\left (  \tau \bu^{n+1/2}  \cdot \nabla \psi,  \frac{\phi^{n+1} - \phi^{n} }{\Delta t} + \bu^{n+1/2}  \cdot \nabla \phi^{n+1/2} \right ) 
+\lambda_1 \left (1, \partial_\phi \rho\, \psi \right ) &= 0, \label{eq:conv_form_ls}  \\
 (1, \rho^{n+1} - \rho^n) &= 0, \label{eq:conv_form_mass_const}
\end{align}
\end{subequations}
where $\mathcal{U} =\{  \mathbf{u} \in  [H^1(\Omega) ]^d ;  \mathbf{u}\cdot\mathbf{n}=0 \} $
and
$ \mathcal{P} =      \{p \in L^2(\Omega); \int p {\rm d}\Omega =0 \} $.
The stabilization parameter and fluid parameters are defined the same way as in 
section \ref{sec:cons_form}. Hence the stability parameter is given by (\ref{eq:tau_def}) and the fluid parameters are given in (\ref{eq:alpha1}) and (\ref{eq:mu_rho_def2}).

\bibliographystyle{unsrt}
\bibliography{references}

\end{document}

%% file: declarations.tex
\def\B#1{\mbox{\boldmath{$#1$}}}

\newcommand{\pd}{\partial}

\newcommand{\nn}{\nonumber}

\newcommand{\bu}{\mathbf{u}}
\newcommand{\bw}{\mathbf{w}}

\newcommand{\bx}{\mathbf{x}}

\newcommand{\onehalf}{\mbox{$\frac{1}{2}$}}

\newcommand{\WW}{\mathcal{W}}

\def\be{\begin{equation}}
\def\ee{\end{equation}}
\def\ba{\begin{array}}
\def\ea{\end{array}}
\def\bea{\begin{eqnarray}}
\def\eea{\end{eqnarray}}
\def\beas{\begin{eqnarray*}}
\def\eeas{\end{eqnarray*}}
\newcommand{\bseq}{\begin{subequations}}
\newcommand{\eseq}{\end{subequations}}
